\documentclass[12pt,letterpaper]{article}
\pdfoutput=1 

\usepackage{etoolbox}
\newtoggle{SUPPLEMENTAL}\toggletrue{SUPPLEMENTAL}
\togglefalse{SUPPLEMENTAL} 

\newtoggle{BLINDED}\toggletrue{BLINDED}
\togglefalse{BLINDED} 
\newcommand{\Blinded}[2]{\iftoggle{BLINDED}{#1}{#2}}

\usepackage[margin=1in,letterpaper]{geometry}
\usepackage{graphicx,grffile}
\usepackage{amsmath,amssymb,amsthm}
\usepackage{mathtools,dsfont}
\usepackage{caption}
\usepackage{subfigure}
\usepackage{booktabs,tabularx,threeparttable,longtable}
\usepackage{siunitx}
\usepackage[shortlabels]{enumitem}
\usepackage[longnamesfirst]{natbib}
\usepackage{alltt}
\usepackage{hypernat}
\usepackage[nodisplayskipstretch]{setspace}
\usepackage{hyperref}
\usepackage[bottom]{footmisc}
\usepackage{xr} 
\usepackage[capitalize,noabbrev]{cleveref}

  \renewenvironment{thebibliography}[1]%
  {\begin{oldthebibliography}{#1}\setlength{\parskip}{0ex}\setlength{\itemsep}{0ex}}%
  {\end{oldthebibliography}}
\appto\TPTnoteSettings{\linespread{1}\footnotesize}

\DeclareGraphicsExtensions{.pdf,.png}
\graphicspath{ {./Figures/} }

\newcommand{\citeposs}[1]{\citeauthor{#1}'s (\citeyear{#1})}
\newcommand{\Citeposs}[1]{\Citeauthor{#1}'s (\citeyear{#1})}

\crefname{conjecture}{Conjecture}{Conjectures}
\crefname{section}{Section}{Sections}
\crefname{subsection}{Section}{Sections}
\crefname{subsubsection}{Section}{Sections}
\Crefname{conjecture}{Conjecture}{Conjectures}
\Crefname{section}{Section}{Sections}
\Crefname{subsection}{Section}{Sections}
\Crefname{subsubsection}{Section}{Sections}
\crefname{appendix}{Appendix}{Appendices}
\crefname{subappendix}{Appendix}{Appendices}
\crefname{subsubappendix}{Appendix}{Appendices}
\Crefname{appendix}{Appendix}{Appendices}
\Crefname{subappendix}{Appendix}{Appendices}
\Crefname{subsubappendix}{Appendix}{Appendices}
\crefname{equation}{}{}
\Crefname{equation}{Equation}{Equations}
\crefname{enumi}{}{}
\Crefname{enumi}{}{}

\crefname{assumption}{}{}
\Crefname{assumption}{Assumption}{Assumptions}

\Crefname{method}{Method}{Methods}

\theoremstyle{plain}
\newtheorem{theorem}{Theorem}

\newtheorem{proposition}[theorem]{Proposition}
\newtheorem{corollary}[theorem]{Corollary}

\theoremstyle{definition}

\newtheorem{assumption}{Assumption}

\newtheorem{method}{Method}

\newcommand{\numnornd}[1]{\num[round-mode=none,group-digits=integer]{#1}} 
\newcommand{\ubar}[1]{\mkern3mu\underline{\mkern-3mu #1\mkern-3mu}\mkern3mu}
\newcommand{\matf}[1]{\ubar{\boldsymbol{\mathbf{#1}}}} 
\newcommand{\vecf}[1]{\boldsymbol{\mathbf{#1}}} 

\newcommand{\SD}[1]{\mathrel{\mathrm{SD}_{#1}}}

\newcommand{\iid}{\stackrel{\mathit{iid}}{\sim}}
\newcommand{\pconv}{\xrightarrow{p}}
\newcommand{\dconv}{\xrightarrow{d}}

\DeclareMathOperator{\Cov}{Cov}
\DeclareMathOperator{\Corr}{Corr}

\newcommand{\R}{{\mathbb R}}

\DeclareMathOperator{\E}{E}
\let\Pr\relax \DeclareMathOperator{\Pr}{P} 

\DeclareMathOperator{\1}{\mathds{1}}
\newcommand{\Ind}[1]{\1\{#1\}}
\newcommand{\NormDist}{\mathrm{N}}
\newcommand{\NormDistp}[2]{\NormDist\left(#1,#2\right)}

\newcommand{\DirDist}{\textrm{Dir}}
\newcommand{\UnifDist}{\textrm{Unif}}
\newcommand{\weaklyto}{\rightsquigarrow}
\newcommand{\diff}{d} 
\newcommand{\pD}[2]{\frac{\partial #1}{\partial #2}}

\newcommand{\closure}[2][3]{%
  {}\mkern#1mu\overline{\mkern-#1mu#2}}
  
\let\originalleft\left
\let\originalright\right
\renewcommand{\left}{\mathopen{}\mathclose\bgroup\originalleft}
\renewcommand{\right}{\aftergroup\egroup\originalright}

\newcommand{\mockalph}[1]{}  
\allowdisplaybreaks[3]

\title{Frequentist properties of {Bayesian} inequality tests%
\thanks{Accepted manuscript for \textit{Journal of Econometrics}, \copyright\ 2020 by the authors. This manuscript version is made available under the CC-BY-NC-ND 4.0 license: \url{http://creativecommons.org/licenses/by-nc-nd/4.0/}}  }

\author{\Blinded{[BLINDED]}{David M.\ Kaplan\thanks{Corresponding author.  Email: \texttt{kaplandm@missouri.edu}. Mail: Department of Economics, University of Missouri, 118 Professional Bldg, 909 University Ave, Columbia, MO 65211-6040, United States. %
Many thanks to Tim Armstrong, Jim Berger, Jeremy Fox, Patrik Guggenberger, 
Jia Li, 
Oliver Linton, Matt Masten, Zack Miller, Stephen Montgomery-Smith, Andriy Norets, Iosif Pinelis, Andres Santos, and anonymous reviewers 
for helpful discussion, comments, examples, and references.
} \qquad Longhao Zhuo\thanks{University of Missouri. Email: \texttt{longhao.zhuo@gmail.com}.}}}%

\date{May 7, 2020} 

\begin{document}

\maketitle

\begin{abstract}
Bayesian and frequentist criteria fundamentally differ, but often posterior and sampling distributions agree asymptotically (e.g., Gaussian with same covariance).
For the corresponding single-draw experiment, we characterize the frequentist size of a certain Bayesian hypothesis test of (possibly nonlinear) inequalities.
If the null hypothesis is that the (possibly infinite-dimensional) parameter lies in a certain half-space, then the Bayesian test's size is $\alpha$; 
if the null hypothesis is a subset of a half-space, then size is above $\alpha$; 
and in other cases, size may be above, below, or equal to $\alpha$.
Rejection probabilities at certain points in the parameter space are also characterized.
Two examples illustrate our results: translog cost function curvature and ordinal distribution relationships.


\end{abstract}

\onehalfspacing

\section{Introduction}
\label{sec:intro}

\subsection{Motivation}

Although Bayesian and frequentist properties fundamentally differ, often we can (approximately) achieve both.
In other cases, Bayesian and frequentist summaries of the data differ greatly.
We characterize the role of null hypothesis ``shape'' in determining such differences.

This paper came out of an empirical question, inspired by \citet{DeatonPaxson1998b,DeatonPaxson1998c}: how can health inequality be assessed when only an ordinal health measure is observed?
Specifically, the ordinal, five-category ``self-reported health status'' is found in many health and economics datasets, but interest is usually in the underlying (latent) health.
Various relationships between two latent distributions can be ``identified'' by various relationships between the corresponding ordinal distributions.
Although ordinal first-order stochastic dominance is relatively simple (i.e., all CDF differences are non-positive), other relationships correspond to more complex subsets of the parameter space (of ordinal CDFs).
Despite such complexity, it is simple to compute Bayesian posterior probabilities for all types of relationships, even simultaneously.
Indeed, often such Bayesian probabilities are straightforward to compute, whereas deriving frequentist tests for general nonlinear inequalities is notoriously difficult \citep[e.g.,][]{Wolak1991}.
This raises the question: do these posterior probabilities have any frequentist interpretation, like a $p$-value?
Existing results comparing Bayesian and frequentist inference do not cover these more complex parameter sets; see \cref{sec:lit}.
For details of this ordinal health example, see \cref{sec:ex-ordinal} and \citet{KaplanZhuo2019b}.

The ordinal health example is a special case of testing whether a certain property holds for all possible partially identified parameter values (or, at least one) within the identified set, as discussed more generally by \citet[p.\ 330]{KlineTamer2016}.
The ordinal distributions can be seen as point-identified parameters that map to identified sets of latent distributions.
Certain ordinal distribution pairs imply that all latent distribution pairs possess some relationship (or, that at least one does).
For example, under certain conditions, latent first-order stochastic dominance implies ordinal first-order stochastic dominance.
Thus, ordinal non-dominance implies latent non-dominance.
The hypothesis of ordinal non-dominance can be interpreted as, ``There is no latent distribution pair in the identified set that satisfies first-order stochastic dominance,'' and ordinal dominance can be interpreted as, ``There is at least one latent distribution pair in the identified set that satisfies first-order stochastic dominance.''

The hypotheses considered in this paper include all the types of hypotheses about the identified set described in \citet[p.\ 330]{KlineTamer2016}.
In addition to the hypotheses in the prior paragraph, they mention the hypothesis that a specified value is in the identified set, which translates to the point-identified parameter being in some region.
They also mention specification testing, where the null hypothesis is that the identified set is nonempty.
This translates to yet another region of the point-identified parameter space, and its complement (where the identified set is empty) is the null hypothesis of \emph{mis}specification.
If the mapping from point-identified parameter to identified set is complex, these regions generally have complex shapes not covered by the existing literature comparing Bayesian and frequentist inference; see \cref{sec:lit}.

Inequality restrictions (often nonlinear) also feature prominently in economic theory.%
\footnote{Nonlinear inequalities come from other sources, too.
          For example, as in \citet{Kaplan2015b}, $H_0 \colon \theta_1\theta_2\ge0$ can test stability of the sign of a parameter over time (or geography), the slope coefficient's sign in IV regression, or whether a treatment attenuates the effect of another regressor.} 
For example, inequalities involving cumulative distribution functions characterize first-order and second-order stochastic dominance, important concepts in welfare analysis and financial portfolio comparisons.
Another example of nonlinear inequalities is curvature constraints on production, cost, indirect utility, and other functions. 
Such constraints usually result from optimization, like utility or profit maximization. 
Some examples are studied in detail in \cref{sec:ex,sec:app-SD1}.
Additional economic examples like bifurcation in dynamic macroeconomic models are reviewed in \cref{sec:results-discussion}.


If (as will be shown) Bayesian and frequentist tests of the null hypothesis can disagree while credible and confidence sets coincide, why not simply report the credible or confidence set?%
\footnote{%
  \Citet{Berger2003} also notes this possible disagreement/agreement, but he writes, ``The disagreement occurs primarily when testing a `precise' hypothesis'' (p.\ 2), whereas we find disagreements even with inequality hypotheses.
  Also, \citet[p.\ 344]{CasellaBerger1987b} opine, ``Interval estimation is, in our opinion, superior to point null hypothesis testing,'' although they do not mention composite null hypotheses like in this paper.}
This makes sense if interest is only in the parameter values themselves.
However, sometimes interest is in testing implications of economic theory or in specification testing.
Other times, inequalities provide economically relevant summaries of a high-dimensional parameter or its properties.

\subsection{Contributions}

Theoretically, we consider a single-draw experiment with equivalent posterior and sampling distributions.
This equivalence often holds asymptotically, as in the various Bernstein--von Mises theorems in the literature.
For example, after centering and scaling, a finite-dimensional parameter estimator often has an asymptotic Gaussian distribution with the same covariance matrix as the asymptotically Gaussian posterior for the same parameter.
We also allow certain non-Gaussian distributions and infinite-dimensional parameters.

To quantify Bayesian--frequentist differences, we characterize the frequentist size and rejection probabilities of a particular Bayesian test.
This test rejects the null hypothesis $H_0$ when the posterior probability of $H_0$ is below $\alpha$, i.e., it treats the posterior like a $p$-value.
Besides being intuitive and practically salient, there are decision-theoretic reasons to examine this test; see \cref{sec:setup-test}.
Although size gives no insight into admissibility, it captures a practical difference between reporting Bayesian and frequentist inferences.

Specifically, we describe how the Bayesian test's size and certain rejection probabilities depend on the shape of the parameter subspace $\Theta_0$ where the null hypothesis is satisfied.
If $\Theta_0$ is a half-space, then the Bayesian test has size $\alpha$, with rejection probability $\alpha$ at every point on the boundary of $\Theta_0$.
If $\Theta_0$ is strictly smaller than a half-space (in a certain sense), then the Bayesian test's size is strictly above $\alpha$, with rejection probability strictly above $\alpha$ at certain boundary points of $\Theta_0$ (like support points).
If $\Theta_0$ is not contained within any half-space, then the Bayesian test's size may be above, equal to, or below $\alpha$, with rejection probability strictly below $\alpha$ at certain boundary points of $\Theta_0$ (like support points of its complement).

\paragraph{Paper structure and notation}

\Cref{sec:lit} reviews the literature.
\Cref{sec:setup} presents the setup and assumptions.
\Cref{sec:results} contains our main results and discussion.
\Cref{sec:ex} illustrates our results with economic examples.
\Cref{sec:app-pfs} contains proofs.
\Cref{sec:app-BvM,sec:app-minimax,sec:app-LFP,sec:app-bivariate,sec:app-translog,sec:app-SD1} contain additional details and examples.
%
Acronyms used include those for 
cumulative distribution function (CDF), 
data generating process (DGP), 
negative semidefinite (NSD), 
posterior expected loss (PEL), 
probability density function (PDF), 
rejection probability (RP), 
and 
single crossing (SC).
%
Notationally, 
$\subseteq$ is subset and $\subset$ is proper subset; 
$\supseteq$ is superset and $\supset$ is proper superset; 
$\mathcal{A}^\complement$ is the complement of set $\mathcal{A}$, and $\closure{\mathcal{A}}$ the closure; 
scalars, (column) vectors, and matrices are respectively formatted as $X$, $\vecf{X}$, and $\matf{X}$; 
$0(\cdot)$ denotes the zero function, i.e., $0(t)=0$ for all $t$.

\section{Literature} 
\label{sec:lit}

Many papers compare Bayesian and frequentist inference.
Here, we highlight examples of different types of conclusions (not all directly comparable to inequality testing): sometimes frequentist inference is more conservative, sometimes Bayesian, sometimes neither.
Then, we discuss three more closely related papers in more detail.

\subsection{Bayesian--frequentist comparisons}
\label{sec:lit-comp}

Some of the literature documents cases where frequentist inference is ``too conservative'' from a Bayesian perspective.
For testing linear inequality constraints of the form $H_0 \colon \vecf{\theta}\ge\vecf{0}$ with $\vecf{\theta}\in\R^d$, $d>1$, \citet{Kline2011} shows frequentist testing to be more conservative, especially with large $d$; see \cref{sec:lit-close} for details.
As another example, under set identification, asymptotically, frequentist confidence sets for the true parameter \citep[e.g.,][]{ImbensManski2004,Stoye2009} are strictly larger than the estimated identified set, whereas Bayesian credible sets are strictly smaller when the prior is over the structural parameter of interest, as shown by \citet[Cor.\ 1]{MoonSchorfheide2012}.\footnote{There seems to be a typo in the statement of Corollary 1(ii), switching the frequentist and Bayesian sets from their correct places seen in the Supplemental Material proof.} 
For testing the null of a unit root in autoregression, \citet{SimsUhlig1991} say frequentist tests ``accept the null more easily'' (p.\ 1592).

Other papers document cases where frequentist inference is ``too aggressive'' from a Bayesian perspective.
Perhaps most famously, in \citeposs{Lindley1957} paradox, the frequentist test can reject when the Bayesian test does not.
\Citet{BergerSellke1987} make a similar argument.
In both cases, as noted by \citet{CasellaBerger1987b}, the results follow primarily from having a large prior probability on a point (or ``small interval'') null hypothesis, specifically $\Pr(H_0)=1/2$.
Arguing that $\Pr(H_0)=1/2$ is ``objective,'' \citet[p.\ 113]{BergerSellke1987} consider even $\Pr(H_0)=0.15$ to be ``blatant bias toward $H_1$.''
\Citet{CasellaBerger1987b} disagree, saying $\Pr(H_0)=1/2$ is ``much larger than is reasonable for most problems'' (p.\ 344).

In yet other cases, Bayesian and frequentist inferences can be similar or even identical, depending on the prior.
For unit root testing, \citet{SimsUhlig1991} determine sample-dependent priors that equate $p$-values and posterior probabilities.
\Citet{CasellaBerger1987} compare Bayesian and frequentist one-sided testing of a location parameter, given a single draw of $X$ from an otherwise fully known density.
They compare the $p$-value to the infimum of the posterior $\Pr(H_0\mid X)$ over various classes of priors.
In many cases, the infimum is attained by the improper prior of Lebesgue measure on $(-\infty,\infty)$ and equals the $p$-value (p.\ 109).
\Citet{BergerEtAl1994} establish an equivalence of Bayesian and \emph{conditional} frequentist testing of a simple null hypothesis against a simple alternative.
\Citet{GoutisEtAl1996} consider jointly testing multiple one-sided hypotheses in a single-draw Gaussian shift experiment.
They propose ``an impartial Bayes rule'' (Section 3.2) that sets the prior $\Pr(H_0)=1/2$ to avoid ``undesirable bias against $H_0$'' (p.\ 1270).
(This also makes the prior odds $\Pr(H_0)/\Pr(H_1)=1$, which simplifies the Bayes factor to the posterior odds, although they only examine the posterior probability of $H_0$.)
If the components of $\vecf{X}$ are mutually independent, then taking the (improper) limit as a symmetric prior's variance goes to infinity, the posterior is proportional to one of the frequentist $p$-values they consider, but it is (weakly) smaller (p.\ 1271).
This complements our setting where we do not impose 
independence, 
$\Pr(H_0)=1/2$, 
finite dimensionality, 
or restricted null hypothesis shape.

There are also many papers showing asymptotic equivalence of Bayesian credible intervals and frequentist confidence sets, due to asymptotic equivalence of posterior and sampling distributions.
Our setup has such an equivalence.
Besides results for point-identified parameters, asymptotically exact frequentist coverage has been established for certain Bayesian credible sets of identified sets, when the prior is (only) over some point-identified parameter and not the partially identified parameter; e.g., see \citet{KlineTamer2016}.
However, these are less directly related to our results on inequality testing.

\subsection{More closely related papers}
\label{sec:lit-close}

\Citet{KlineTamer2016} mention many of the issues related to our paper, although they do not focus on frequentist properties of Bayesian tests.
Their Sections 3 and 4 describe and interpret different statements about the identified set, but the focus is on computation and consistency of Bayesian posterior probabilities.
Asymptotic approximations based on normal variables are given (Theorems 1(iii) and 3(iii,vi)) but not compared to frequentist $p$-values.
In Section 5, in the context of a certain credible set, they raise the question of equality of posterior probability and $p$-value (pp.\ 348--349).
They state this holds for interval-identified parameters but otherwise just refer to \citet{Kline2011}.

\Citet{Kline2011} studies a special case of our setting.
Similar to our setting (although we allow certain non-normal distributions), his posterior is multivariate normal (Definition 1, p.\ 3133), with covariance matrix equal to that of the multivariate normal sampling distribution.
\Citet{Kline2011} motivates this equivalence from estimating a vector of means, where either the likelihood is normal and the prior is flat, or there is a Dirichlet process prior under which the asymptotic posterior is equivalent to the asymptotic sampling distribution (Sections 2.1--2.2).
Most restrictive compared to our setting, the null hypothesis is that all components of the vector are non-negative.

Most relevant to our paper are \citeposs{Kline2011} Lemmas 1 and 4 (pp.\ 3134, 3140).
They show that if the null hypothesis is satisfied and at least two components of the mean vector equal zero (on the boundary of the null), then the asymptotic sampling distribution of the posterior probability of the null is first-order stochastically dominated by the standard uniform distribution.
That is, if $\hat\pi$ is the posterior probability of $H_0$, then $\Pr(\hat\pi \le \alpha)>\alpha$ for all $\alpha\in(0,1)$.
Consequently, the test that rejects when the posterior is below $\alpha$ has type I error rate (and thus size) above $\alpha$ in such a situation.
%
Although not discussed in the original paper, this also implies $\Pr(\hat\pi \ge 1-\alpha)<\alpha$, so the type I error rate is below $\alpha$ if the null and alternative hypotheses are switched.
\Citet[\S4]{Kline2011} also provides numerical examples and simulations showing that the magnitude of the discrepancy increases with the dimension.
%

\Citeposs{Kline2011} results are consistent with ours, but we study a more general null hypothesis.
By restricting the null hypothesis, \citet{Kline2011} can more precisely quantify the Bayesian test's rejection probabilities (Lemma 1) and can propose a frequentist test using the posterior probability of $H_0$ as the test statistic, but with critical value smaller than $\alpha$ (which can be interpreted as an adjusted loss function, in light of our \cref{sec:setup-test}).
Our null hypothesis is simply that the parameter lies in some subset of the parameter space, with no restrictions on the shape of the subset.
Consequently, our results apply to more settings.
For example, the applications in \cref{sec:ex} and most others mentioned in \cref{sec:intro,sec:results} are not covered by the results of \citet{Kline2011}.
In particular, we allow hypotheses about identified sets without restricting the mapping from point-identified parameter to identified set.

\Citet{EfronTibshirani1998} consider the confidence level for a parameter $\vecf{\mu}$ being in a certain region of the parameter space.
They call this the ``problem of regions.''
In particular, they focus on the confidence level of the region in which the point estimate $\hat{\vecf{\mu}}$ lies.
They note (p.\ 1698) that their ``first-order bootstrap'' confidence value is essentially a Bayesian posterior probability given a flat prior.
Although their main contribution is a bootstrap adjustment to better match the ``correct'' frequentist confidence level, most related to our paper is their discussion of how the Bayesian (i.e., first-order bootstrap) and frequentist measures differ depending on the curvature of the boundary between the regions.
Qualitatively, they note the Bayesian posterior probability is higher than the frequentist confidence level when the boundary ``curves away from $\hat{\vecf{\mu}}$,'' whereas it is lower when the boundary ``curves toward $\hat{\vecf{\mu}}$'' (p.\ 1690).
In terms of hypothesis testing, this implies the Bayesian test is more likely (than a frequentist test) to reject the null when the boundary of the null region ``curves away from $\hat{\vecf{\mu}}$,'' i.e., has a convex boundary near $\hat{\vecf{\mu}}$.
Similarly, it implies the Bayesian test is less likely to reject if the null region's boundary is non-convex near $\hat{\vecf{\mu}}$.
\Citet{EfronTibshirani1998} also discuss the better agreement of Bayesian and frequentist confidence when i) the boundary can be written as a level set of a smooth function of the underlying parameter and ii) the Welch--Peers prior \citep{WelchPeers1963} can be used (e.g., \S4 or Remark B on p.\ 1714).

However, the theoretical results in \citet{EfronTibshirani1998} stipulate an asymptotic approximation local to a smooth boundary of the null region.
More precisely, up to smaller-order terms, the boundary is treated as quadratic, as in their (2.10).
This greatly simplifies matters; e.g., the distance from $\hat{\vecf{\mu}}$ to the nearest boundary point has a normal distribution with mean and variance related to the quadratic coefficient matrix, as in their (2.12).
Asymptotically, the smooth boundary looks nearly flat, with only $O(n^{-1/2})$ curvature (p.\ 1692 before (2.14)).
Consequently, even with no adjustment, Bayesian and frequentist confidence levels differ by merely $O_p(n^{-1/2})$; see (2.24) and following text (p.\ 1694).
Further, the approximate connection with Welch--Peers priors in Section 4 assumes the boundary is a level set of some scalar-valued function of $\vecf{\mu}$.
In contrast, our results do not restrict the boundary of the null hypothesis region in any way.
This is important since many null hypotheses entail non-smooth boundaries, such as joint tests of multiple inequalities \citep[as in][]{Kline2011}.
Of course, the benefit of their stronger assumptions is that \citet{EfronTibshirani1998} can more precisely approximate the quantitative difference between frequentist confidence and Bayesian posterior probability, as in their (2.21) or (2.24).

\section{Setup and assumptions}
\label{sec:setup}

In \cref{sec:setup-test}, a specific Bayesian hypothesis test is described along with the decision-theoretic context. 
In \cref{sec:setup-assumptions}, the assumptions used for the results in \cref{sec:results} are presented and discussed.
\Cref{sec:setup-LFP} discusses an alternative approach for this setting.

\subsection{The {Bayesian} test and decision-theoretic context}
\label{sec:setup-test}

The Bayesian test rejects the null hypothesis when its posterior probability is below $\alpha$.
\begin{method}[Bayesian test]\label{meth:Bayes}
Reject $H_0$ if $\Pr( H_0 \mid \textrm{data}) \le \alpha$; otherwise, accept $H_0$. 
\end{method}

Besides seeming intuitive, \cref{meth:Bayes} is a generalized Bayes decision rule.
Specifically, it minimizes posterior expected loss (PEL) for the loss function taking value $1-\alpha$ for type I error, $\alpha$ for type II error, and zero otherwise.
Letting $H=0$ if $H_0$ is true and $H=1$ if $H_0$ is false, the loss function is
\begin{equation}\label{eqn:L}
\begin{split}
L(H=0, \textrm{reject}) &= 1-\alpha , \quad
L(H=1, \textrm{accept}) = \alpha , \\
L(H=1, \textrm{reject}) &= 
L(H=0, \textrm{accept}) = 0 .
\end{split}
\end{equation}
Thus, given a particular decision but unknown $H$,
\begin{equation}\label{eqn:L-decision}
\begin{split}
L(H,\textrm{reject})
&= (1-\alpha) \Ind{H=0}
  +(0)\Ind{H=1}
 = (1-\alpha) \Ind{H_0\textrm{ true}} , \\
L(H,\textrm{accept})
&= (0) \Ind{H=0}
  +(\alpha)\Ind{H=1}
 = \alpha \Ind{H_0\textrm{ false}} .
\end{split}
\end{equation}
To compute PEL, let $\Pr(\cdot \mid \vecf{X})$ denote the posterior probability given observed data $\vecf{X}$, and $\E(\cdot\mid\vecf{X})$ the posterior expectation.
The PEL for the decision to reject $H_0$ is
\begin{equation}\label{eqn:PEL-rej}
\E[L(H, \textrm{reject}) \mid \vecf{X}]
= \E[\overbrace{(1-\alpha)\Ind{H_0\textrm{ true}}}^{\textrm{from \cref{eqn:L-decision}}}
    \mid \vecf{X}]
= (1-\alpha)\Pr(H_0 \mid \vecf{X}) ,
\end{equation}
i.e., the type I error loss ($1-\alpha$) times the posterior probability that $H_0$ is true and thus ``reject'' is a type I error.
Similarly, the PEL of accepting $H_0$ is
\begin{equation}\label{eqn:PEL-acc}
\E[L(H, \textrm{accept}) \mid \vecf{X}]
= \E[\overbrace{\alpha\Ind{H_0\textrm{ false}}}^{\textrm{from \cref{eqn:L-decision}}} \mid \vecf{X}]
= \alpha[1-\Pr(H_0\mid \vecf{X})],
\end{equation}
i.e., the type II error loss ($\alpha$) times the probability that $H_0$ is false and thus ``accept'' is a type II error.
PEL is thus minimized by rejecting $H_0$ if $\Pr(H_0 \mid \vecf{X}) \le \alpha$ and accepting $H_0$ otherwise; i.e., \cref{eqn:PEL-rej} is smaller than \cref{eqn:PEL-acc} if and only if $\Pr(H_0 \mid \vecf{X})\le\alpha$, with equality when $\Pr(H_0 \mid \vecf{X})=\alpha$.

The Bayesian test's type I error rates and size are compared to $\alpha$ (instead of some other value) for both practical and decision-theoretic reasons.
Practically, the Bayesian test can be seen as treating the posterior as a (frequentist) $p$-value; we want to know if this is valid, similar to \citet{CasellaBerger1987}. 
Decision-theoretically, given the same loss function from \cref{eqn:L} 
under which the Bayesian test in \cref{meth:Bayes} is a generalized Bayes rule, any unbiased frequentist test with size $\alpha$ is a minimax risk decision rule; see \cref{sec:app-minimax}.

Ideally, a single decision rule minimizes both maximum risk and PEL.
However, if the Bayesian test's size is above or below $\alpha$, this is not possible.

\subsection{Assumptions}
\label{sec:setup-assumptions}

The two main assumptions are stated and discussed in \cref{sec:a-multi,sec:a-H0}.

\subsubsection{Posterior and sampling distributions}
\label{sec:a-multi}

\Cref{a:multi} states conditions on the sampling and posterior distributions of a functional $\phi(\cdot)$ in the single-draw experiment we consider. 
As usual, the sampling distribution treats the (functional of the) data $\phi(\vecf{X})$ as random and conditions on the parameter $\vecf{\theta}$, whereas the posterior treats the (functional of the) parameter $\phi(\vecf{\theta})$ as random and conditions on the data $\vecf{X}$.

\begin{assumption}\label{a:multi}
Let $F(\cdot)$ be a continuous CDF with support $\R$ and symmetry $F(-x)=1-F(x)$. 
Let the single observation $\vecf{X}$ and the parameter $\vecf{\theta}$ both belong to a Banach space of possibly infinite dimension. 
Let $\phi(\cdot)$ denote a continuous linear functional, with sampling distribution 
$\phi(\vecf{X})-\phi(\vecf{\theta}) \mid \vecf{\theta} \sim F$ 
and posterior 
$\phi(\vecf{\theta})-\phi(\vecf{X})\mid \vecf{X}\sim F$, for every $\vecf{\theta}$ and every $\vecf{X}$ in the space.
\end{assumption}

The hope is that the conditions in \cref{a:multi} are reasonable approximations of the conditions encountered in practice.
The use of a single draw of $\vecf{X}$ has precedent in papers like \citet{CasellaBerger1987} and \citet{GoutisEtAl1996}.
Usually, $\vecf{X}$ in \cref{a:multi} represents some estimator $\hat{\vecf{\theta}}$.
For example, in a particular empirical setting, if the estimator (centered at the true parameter) and posterior (centered at the estimate) are both approximately Gaussian with approximately the same covariance, then our results should provide some insight.
The conditions in \cref{a:multi} are now further discussed.

\Cref{a:multi} includes the common case of a Gaussian distribution.
In $\R^d$, continuous linear functionals are simply linear combinations $\phi(\vecf{X}) = \vecf{c}'\vecf{X}$ for constant vector $\vecf{c}\in\R^d$.
If $\vecf{X}\sim\NormDist(\vecf{\theta},\matf{V})$ and $\vecf{\theta}\sim\NormDist(\vecf{X},\matf{V})$, then the sampling distribution of $\vecf{c}'\vecf{X}-\vecf{c}'\vecf{\theta}$ and the posterior distribution of $\vecf{c}'\vecf{\theta}-\vecf{c}'\vecf{X}$ are both $\NormDist(0,\vecf{c}'\matf{V}\vecf{c})$, satisfying the equivalence and symmetry conditions in \cref{a:multi}.
In infinite-dimensional spaces, if $X(\cdot)$ is a Gaussian process in some Banach space and $\phi(\cdot)$ belongs to the dual of that space, then $\phi\bigl(X(\cdot)\bigr)$ is a scalar Gaussian random variable; e.g., see Definition 2.2.1(ii) in \citet[p.\ 42]{Bogachev1998} and \citet[pp.\ 376--377]{vanderVaartWellner1996}.
For example, given Gaussian process $X(\cdot)$ and $r,r_1,r_2\in\R$, both $\phi(X(\cdot))=X(r)$ and $\phi(X(\cdot))=X(r_2)-X(r_1)$ are Gaussian (relevant for first-order stochastic dominance and monotonicity, respectively).

The equivalence of posterior and sampling distributions in \cref{a:multi} requires that the prior not influence the posterior.
Implicitly, \cref{a:multi} entails an improper uninformative (flat) prior.
For example, take scalar $X,\theta\in\R$, so $\phi(X)=X$.
With sampling distribution $X \mid \theta \sim \NormDistp{\theta}{1}$ and prior $\theta\sim\NormDist(m,\tau^2)$, the posterior is 
\begin{equation*}
\theta \mid X 
\sim \NormDistp{\frac{\tau^2 X+m}{\tau^2+1}}{\frac{\tau^2}{\tau^2+1}} .
\end{equation*}
The equivalence condition is satisfied (only) in the limit as $\tau^2\to\infty$, yielding the posterior $\theta\mid X \sim \NormDistp{X}{1}$.
Alternatively, the improper ``flat'' prior (equal to $1$ for all $\theta$) yields the same posterior.
This happens because the posterior is proportional to the likelihood times the prior, so the flat prior disappears and leaves the posterior proportional to the likelihood, which in this example is $(2\pi)^{-1/2}\exp\{-(X-\theta)^2/2\}$.
Interpreted as a posterior, this gives the PDF of $\theta$ given $X$ as that of a $\NormDist(X,1)$ distribution.
This improper prior is fine for \cref{meth:Bayes} because only posterior probabilities are used, unlike with Bayes factors that involve a ratio of prior probabilities \citep[e.g.,][]{BayarriEtAl2012}.

The uninformative prior implicit in \cref{a:multi} is restrictive but often relevant.
It is restrictive because in practice the prior may influence the posterior greatly.
Even in large samples, a sufficiently strong prior can always influence the posterior.
However, to gain any general understanding, there must be some restriction on the prior, otherwise it can arbitrarily affect the results (e.g., a dogmatic prior on a point inside $\Theta_0$, or outside).
Further, in practice, often the prior does not influence the posterior much.
This may be due to a large sample and/or an uninformative prior on $\vecf{\theta}$.
For example, in the application to partially identified parameters in a binary entry game, \citet[p.\ 355]{KlineTamer2016} specify ``an uninformative conjugate Dirichlet prior,'' similar to that in our application in \cref{sec:ex-ordinal}.
Finally, our results can be seen as trying to isolate the effect of the Bayesian or frequentist framework, apart from any particular prior.

\Cref{a:multi} can be interpreted asymptotically.
The asymptotic equivalence of posterior and sampling distributions is the conclusion of a Bernstein--von Mises theorem.
Many Bernstein--von Mises theorems have been established in different parametric, semiparametric, and even nonparametric settings; see \cref{sec:app-BvM}.
Most directly, such results could approximately justify the posterior and sampling distribution equivalence in \cref{a:multi}, where $\vecf{\theta}$ is the parameter of interest estimated by $\hat{\vecf{\theta}}=\vecf{X}$.
Alternatively, if \cref{a:multi} describes a limit experiment, then $\vecf{\theta}$ is a \emph{local} mean parameter. 
Then, the ``size'' of the test is really the ``local size,'' i.e., the supremum type I error in a local neighborhood of some parameter value; but statements about rejection probability retain the same interpretation.
Also, the shape of the null hypothesis region may change substantially when ``zooming in'' on a particular local neighborhood; see below for further discussion.

The following is a simple example of an asymptotic limit experiment with a local parameter.
Let $X,\theta\in\R$ with $\phi(v)=v$.
If $Y_{ni} \iid \NormDist(\mu_n,1)$, $i=1,\ldots,n$, and $\sqrt{n}\mu_n\to\theta$, then $\sqrt{n}\bar{Y}_n=n^{-1/2}\sum_{i=1}^{n}Y_{ni} \dconv \NormDistp{\theta}{1}$.
More generally, if 
$Y_{ni} \iid \NormDistp{m+\mu_n}{\sigma^2}$ 
and 
$\sqrt{n}\mu_n\to\theta$, 
then
\begin{equation}\label{eqn:local-scalar}
X_n \equiv 
\sqrt{n}(\bar{Y}_n-m)/\hat\sigma
\dconv X
\sim \NormDist(\theta,1)
\end{equation}
given consistent estimator $\hat\sigma^2 \pconv \sigma^2$.
Since $\theta$ is the \emph{local} mean parameter, assuming $\theta \in \R$ does not require that $\R$ is the original parameter space (e.g., the space for $m+\mu_n$ in the example), but it does exclude boundary points.
This type of result holds for a wide variety of models, estimators, and sampling assumptions.
It is most commonly used for local power analysis but has been used for purposes like ours in papers like \citet[eqn.\ (4.2)]{AndrewsSoares2010}.

For our purpose of approximating the finite-sample frequentist size of a Bayesian test, considering a fixed DGP and drifting centering parameter can be just as helpful as considering a fixed centering parameter and drifting DGP.
This allows the relevant Bernstein--von Mises theorem to hold only for fixed (not drifting) DGPs.
For example, in $\R^d$,
\begin{equation}\label{eqn:gen-drift-mu0}
\vecf{X}_n 
 = \sqrt{n}(\hat{\vecf{\mu}}-\vecf{\mu}_{0,n}) 
 = \overbrace{\sqrt{n}(\hat{\vecf{\mu}}-\vecf{\mu}) }^{\dconv \NormDistp{\vecf{0}}{\matf{\Sigma}}}
          +\overbrace{\sqrt{n}(\vecf{\mu}-\vecf{\mu}_{0,n})}^{=\vecf{\theta}_n\to\vecf{\theta}}
\dconv 
\overbrace{%
\vecf{X} \sim \NormDistp{\vecf{\theta}}{\matf{\Sigma}}
}^{\textrm{limit experiment}}
,
\end{equation}
where $\vecf{\mu}$ is the true parameter value, $\hat{\vecf{\mu}}$ is a $\sqrt{n}$-normal estimator (as seen), $\matf{\Sigma}$ is the asymptotic covariance matrix (known or consistently estimable), 
and 
$\vecf{X}_n$ is a statistic based on $\hat{\vecf{\mu}}$ and centered at the deterministic sequence $\vecf{\mu}_{0,n}$ that satisfies $\sqrt{n}(\vecf{\mu}-\vecf{\mu}_{0,n})  \to  \vecf{\theta}$, the local mean parameter.
This does not have a literal meaning like ``we must change $\vecf{\mu}_0$ if our sample size increases,'' just as a drifting DGP does not mean literally that ``the population distribution changes as we collect more data''; rather, it is simply a way to capture the idea of $\vecf{\mu}_0$ being ``close to'' the true $\vecf{\mu}$ in the asymptotics.
Analogously, for the posterior, letting 
$\vecf{\theta}_n = \sqrt{n}(\vecf{\mu}-\vecf{\mu}_{0,n})$ 
and again 
$\vecf{X}_n = \sqrt{n}(\hat{\vecf{\mu}}-\vecf{\mu}_{0,n})$, 
assuming a Bernstein--von Mises theorem,
\begin{equation}\label{eqn:gen-drift-mu0-posterior}
\vecf{\theta}_n - \vecf{X}_n
= 
\sqrt{n}(\vecf{\mu}-\hat{\vecf{\mu}}) \dconv \NormDistp{\vecf{0}}{\matf{\Sigma}} 
.
\end{equation}

\subsubsection{Shape of null hypothesis region}
\label{sec:a-H0}

\begin{assumption}\label{a:H0}
The null hypothesis is $H_0 \colon \vecf{\theta} \in \Theta_0$, where $\Theta_0$ is a subspace of the Banach space in \cref{a:multi}; the alternative hypothesis is that $\vecf{\theta} \not\in \Theta_0$.
\end{assumption}

\Cref{a:H0} is very general.
The shape of $\Theta_0$ is not restricted at all.
It does assume each $\vecf{\theta}$ satisfies either the null or alternative hypothesis (i.e., there is no third region), but this is often true.

In practice, it matters whether $\vecf{\theta}$ is interpreted as the original parameter or as a local parameter as in \cref{eqn:local-scalar} or \cref{eqn:gen-drift-mu0-posterior}.
In our results, the shape of $\Theta_0$ matters; should this be the ``global shape'' of the null hypothesis region, or the ``local shape'' in a neighborhood of the true parameter?
The answer is not obvious.
In practice, one idea is to construct something like a $99.9\%$ credible set or confidence set to get an idea of the ``local'' shape, where ``local'' is based on the data, similar in spirit to the idea of \citet{RomanoShaikhWolf2014}.

Consider some examples in which the relevant shape of $\Theta_0$ depends on the sample size as well as the true parameter value.
Let $\vecf{X},\vecf{\theta}\in\R^2$.
First, let $H_0 \colon \theta_1\le0,\theta_2\le0$.
Thus, $\Theta_0$ is the third quadrant, with a convex corner at the origin and flat sides extending down each axis to $-\infty$.
If we look at a neighborhood of the origin (or equivalently scale by $\sqrt{n}$), the shape is identical; no matter how much we zoom in, the null hypothesis is the third quadrant since $(\theta_1,\theta_2)\le(0,0) \iff (\sqrt{n}\theta_1,\sqrt{n}\theta_2)\le(0,0)$.
However, consider a neighborhood of $(\theta_1,\theta_2)=(-m,0)$.
Perhaps with small sample size $n$, the standard deviation is much larger than $m$, so the original $\Theta_0$ shape with the corner is appropriate.
But with large enough $n$, as the standard deviation shrinks to zero, a $n^{-1/2}$-neighborhood around $(-m,0)$ no longer includes the corner at the origin.
Instead, the null hypothesis boundary appears flat.
That is, it becomes obvious that $\theta_1<0$, so the null hypothesis reduces to $H_0 \colon \theta_2 \le 0$.
This is the same idea behind moment inequality testing procedures like in \citet{RomanoShaikhWolf2014}.
They hope to learn which inequalities might be binding, which is equivalent to the local shape of $\Theta_0$.

Consider a variant of the above example.
Remove a notch from $\Theta_0$ near the origin: if $\theta_1>-r$ and $\theta_2>-r$, then the point $(\theta_1,\theta_2)$ is excluded from $\Theta_0$ even if $(\theta_1,\theta_2)\le(0,0)$.
If $r$ is very small compared to the standard deviation, then this notch's effect is negligible; $\Theta_0$ still looks like the third quadrant.
However, if the sample size is large enough that $r$ is much bigger than the standard deviation, a local neighborhood of $(-r,-r)$ looks very different: the null hypothesis is satisfied in quadrants 2, 3, and 4, not just quadrant 3.
Instead of a convex corner, the null region has a non-convex corner, which leads to very different behavior of the Bayesian test.

This difference between global and local shape (and its practical importance) is illustrated further in the example of \cref{sec:ex-ordinal}.

\subsection{Least favorable prior approach}
\label{sec:setup-LFP}

Comparing Bayes and minimax decision rules evokes the least favorable prior.
Most simply, if there exists a (least favorable) prior such that the Bayes risk of the corresponding Bayes rule equals the supremum risk (over the parameter space) of the same Bayes rule, then that rule is minimax \citep[Thm.\ 1.4 in Ch.\ 5, p.\ 310]{LehmannCasella1998}.
More generally, there is a least favorable \emph{sequence} of priors, with a corresponding sequence of Bayes rules and their Bayes risks (against the prior); any decision rule whose supremum risk equals the limit of any such Bayes risk sequence is minimax \citep[Thm.\ 1.12 in Ch.\ 5, p.\ 316]{LehmannCasella1998}.

Deriving a least favorable prior (sequence) would be valuable, although there are two reasons we do not pursue it in depth here.
First, for hypothesis testing, the (limiting) least favorable prior is never the flat prior implicit in \Cref{a:multi}.
Even if such a prior is found, it is still helpful to understand the Bayesian test's behavior under \cref{a:multi}.
Second, for very general problems like in \Cref{a:H0}, minimax results are notoriously difficult to obtain.
For example, before discussing least favorable priors, \citet[\S5.1, p.\ 309]{LehmannCasella1998} remark, ``We shall not be able to determine minimax estimators for large classes of problems but, rather, will treat problems individually.''
(Hypothesis testing is essentially ``estimation'' of the population object $\Ind{\vecf{\theta}\not\in\Theta_0}$.)

For reference, a simple example is worked out in \cref{sec:app-LFP}.
The null hypothesis is $H_0 \colon \theta\le0$, with scalar $X\sim\NormDist(\theta,1)$.
Despite the simple setup, it is more complicated than a simple estimation problem, where any sequence of priors approaching the flat (improper) prior is least favorable.
Unlike with estimation, the testing decision is most difficult near $\theta=0$, so the least favorable prior sequence must (in the limit) concentrate mass near that point.
Otherwise, if for example there is significant prior mass on $\theta>10$, it is obvious that the test should reject in such cases, which lowers the posterior expected loss below the minimax risk.
The prior also cannot put all mass on $\theta=0$, otherwise the Bayesian test never rejects and has zero Bayes risk.
One least favorable prior sequence places mass on only two points: $\theta=0$, and a strictly positive value.
In the limit, the strictly positive value converges to zero, and the prior probability of $\theta=0$ converges to $\alpha$ in a particular way (depending on the value of the positive support point).
We show the limit of the corresponding Bayes risk sequence equals the minimax risk $\alpha(1-\alpha)$ from \cref{sec:app-minimax}.
Details and discussion are in \cref{sec:app-LFP}.

Although the example in \cref{sec:app-LFP} works out intuitively, it is unclear how well this approach generalizes.
Nonetheless, it is worth investigation in future work.

\section{Results and discussion}
\label{sec:results}

\Cref{thm:1} contains this paper's main results. 
Discussion and special cases follow.

\begin{theorem}\label{thm:1}
Let \Cref{a:multi,a:H0} hold.
Consider the Bayesian test of \cref{meth:Bayes}.
\begin{enumerate}[label={\textup{(\roman*)}}, ref=\cref{thm:1}(\roman*)]
\item\label{thm:half} For any $\phi(\cdot)$ satisfying \cref{a:multi}, if there exists $c_0 \in \R$ such that 
$\Theta_0 = \{\vecf{\theta} : \phi(\vecf{\theta}) \le c_0\}$ (i.e., $\Theta_0$ is a half-space), 
then the Bayesian test's size is $\alpha$, and its type I error rate is $\alpha$ when $\phi(\vecf{\theta})=c_0$. 
\item\label{thm:subeq} 
For any $\phi(\cdot)$ satisfying \cref{a:multi}, if there exists $c_0 \in \R$ such that 
$\Theta_0 \subseteq \{\vecf{\theta} : \phi(\vecf{\theta}) \le c_0\}$ 
with $c_0\in\phi(\closure{\Theta_0})$ (where $\closure{\Theta_0}$ denotes the closure), 
then the Bayesian test's size is $\alpha$ or greater. 
\item\label{thm:sub} Continuing from \cref{thm:subeq}, further assume that there exists a $\phi_2(\cdot)$ satisfying $\phi_2(\vecf{X})-\phi_2(\vecf{\theta})\mid\vecf{\theta}\sim F_2$ and $\phi_2(\vecf{\theta})-\phi_2(\vecf{X})\mid\vecf{X}\sim F_2$ for every $\vecf{\theta}$ and $\vecf{X}$ and satisfying the properties below.
Assume the $\phi(\cdot)$ from \cref{thm:subeq} may be written as $\phi(\cdot)=\phi_3(\phi_2(\cdot))$ for some projection $\phi_3(\cdot)\colon\R^d\mapsto\R$.
Assume there exists (in $\R^d$) a set
$\Phi_2  \equiv  \{ \phi_2(\vecf{\theta}) : \vecf{\theta} \in \Theta_0 \}  \subset  \{ \phi_2(\vecf{\theta}) : \phi(\vecf{\theta}) \le c_0 \}$.
Further assume 
(a) the set $\Delta \equiv \{ \vecf{p} : \vecf{p} \in \R^d, \phi_3(\vecf{p}) \le c_0, \vecf{p} \not\in \Phi_2 \}$ has positive Lebesgue measure, 
(b) $F_2$ has a strictly positive PDF over $\R^d$, and 
(c) $\Pr(\phi_2(\vecf{\theta}) \in \Phi_2 \mid \phi_2(\vecf{X}))$ is continuous in $\phi_2(\vecf{X})$. 
Then, the Bayesian test's rejection probability is strictly above $\alpha$ when $\phi(\vecf{\theta})=c_0$, and its size is strictly above $\alpha$.
If $\vecf{\theta} \in \R^d$, then $\phi_2(\vecf{\theta})=\vecf{\theta}$, $\phi_3(\cdot)=\phi(\cdot)$, and $\Phi_2=\Theta_0$, 
and the rejection probability is strictly above $\alpha$ for any $\vecf{\theta}$ that is a support point of the closure of $\Theta_0$.
\item\label{thm:non} If, contrary to \cref{thm:subeq}, there do not exist any $\phi(\cdot)$ and $c_0$ such that $\Theta_0 \subseteq \{\vecf{\theta} : \phi(\vecf{\theta}) \le c_0\}$ (i.e., $\Theta_0$ is not a subset of a half-space), 
then the Bayesian test's size may be greater than, equal to, or less than $\alpha$, and it may depend on the sampling/posterior distribution (in addition to the shape of $\Theta_0$).
\item\label{thm:super} Continue the notation and conditions from \cref{thm:sub}, but redefine $\Phi_2 \equiv \{\phi_2(\vecf{\theta}) : \vecf{\theta}\not\in\Theta_0\}^\complement \supset \{ \phi_2(\vecf{\theta}) : \phi(\vecf{\theta}) \le c_0\}$, again with $\Phi_2=\Theta_0$ if $\vecf{\theta}\in\R^d$; and redefine $\Delta \equiv \{\vecf{p} : \vecf{p}\in\Phi_2, \phi_3(\vecf{p})>c_0 \}$ (the ``extra'' part of $\Phi_2$ beyond the half-space in $\R^d$).
Then, the Bayesian test's rejection probability is strictly below $\alpha$ when $\phi(\vecf{\theta})=c_0$.
With $\vecf{\theta}\in\R^d$, this implies the rejection probability is strictly below $\alpha$ for any $\vecf{\theta}$ that is a support point of the closure of the complement of $\Theta_0$, i.e., a support point of $\closure{\Theta_0^\complement}$.
\end{enumerate}
\end{theorem}

\subsection{Intuition for results}
\label{sec:results-intuition}

Some intuition for \cref{thm:1} is now provided.
The focus is on the finite-dimensional case, specifically $\R^2$; the infinite-dimensional argument essentially finds a finite-dimensional projection in which the same logic goes through.
For simplicity, here the inverse CDF $F^{-1}(\cdot)$ is assumed to exist.

The one-dimensional version of \cref{thm:half} is well known but worth mentioning since it underlies all other results.
The key is the symmetry of $F(\cdot)$ in \Cref{a:multi}.
Without symmetry, as seen below, the Bayesian test's size would not generally equal $\alpha$.
Consider $H_0 \colon \theta\le0$.
The posterior probability of $H_0$ is strictly decreasing in $X$, and it equals $\alpha$ when $X=-F^{-1}(\alpha)$, as shown below, so the Bayesian test rejects when $X>-F^{-1}(\alpha)$.
To verify this critical value, since the posterior is $\theta-X\mid X\sim F$, and $\theta\le0$ can be written $\theta-X\le-X$,
\begin{align*}
\Pr(H_0 \mid X = -F^{-1}(\alpha))
  &= \Pr(\theta\le0 \mid X = -F^{-1}(\alpha))
\\&= \Pr(\theta-X \le -X \mid X = -F^{-1}(\alpha))
\\&= F(-(-F^{-1}(\alpha)))
\\&= F(F^{-1}(\alpha))
\\&= \alpha .
\end{align*}
The power function is increasing in $\theta$, so the type I error rate is largest at $\theta=0$.
Since the sampling distribution of $X-\theta$ (which is just $X$ when $\theta=0$) also has CDF $F(\cdot)$, size is
\begin{align*}
\sup_{\theta\le0}\Pr(X > -F^{-1}(\alpha) \mid \theta)
  &= \Pr(X > -F^{-1}(\alpha) \mid \theta=0)
\\&= 1 - \Pr(X-\theta \le -F^{-1}(\alpha)-\theta \mid \theta=0)
\\&= 1 - F(-F^{-1}(\alpha)) .
\end{align*}
Without symmetry, this does not equal $\alpha$.
With the symmetry condition $F(-x)=1-F(x)$ in \cref{a:multi}, however, the Bayesian test's size simplifies to
\begin{equation*}
1 - F(-F^{-1}(\alpha))
= 1- [1-F(F^{-1}(\alpha))]
= \alpha .
\end{equation*}

\begin{figure}[htb]
\centering
\subfigure[Illustration of \cref{thm:half}]{\label{fig:thm-intuition-half}%
\includegraphics{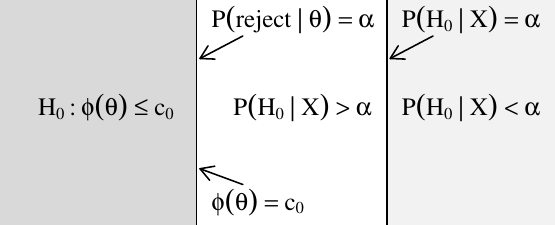}
}
\subfigure[Illustration of \cref{thm:sub}]{\label{fig:thm-intuition-sub}%
\includegraphics{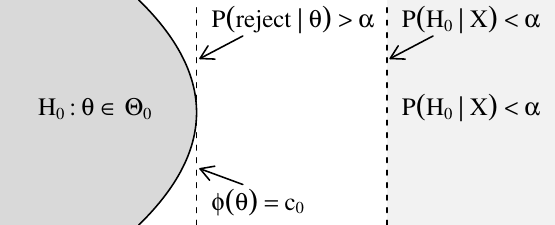}
}
\subfigure[Illustration of \cref{thm:super}]{\label{fig:thm-intuition-super}%
\includegraphics{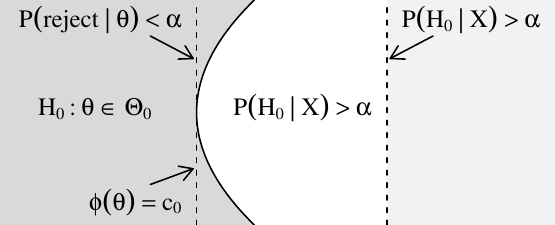}
}
\caption{\label{fig:thm-intuition}Illustrations of \cref{thm:1} for $\vecf{X},\vecf{\theta}\in\R^2$.}
\end{figure}

\Cref{fig:thm-intuition} illustrates the intuition behind \cref{thm:1}.
Here, $\vecf{X},\vecf{\theta}\in\R^2$, so $\phi_2(\vecf{\theta})=\vecf{\theta}$, $\phi_3(\cdot)=\phi(\cdot)$, and $\Phi_2=\Theta_0$.

\Cref{fig:thm-intuition-half} illustrates \cref{thm:half}.
Here, $\vecf{X}=(X_1,X_2)$ and $\vecf{\theta}=(\theta_1,\theta_2)$, where $X_1$ and $\theta_1$ are on the horizontal axis while $X_2$ and $\theta_2$ are on the vertical axis.
Further, $\phi(\vecf{\theta})=\theta_1$, and imagine $c_0=0$ for simplicity.
The left vertical line shows $\phi(\vecf{\theta})=c_0$, i.e., $\theta_1=0$.
Thus, the null is $H_0 \colon \theta_1\le0$, reducing to the one-dimensional case.
\Cref{a:multi} says $\phi(\vecf{X})-\phi(\vecf{\theta})\mid\vecf{\theta}\sim F$ and $\phi(\vecf{\theta})-\phi(\vecf{X})\mid\vecf{X}\sim F$, which here is $X_1-\theta_1\mid\vecf{\theta}\sim F$ and $\theta_1-X_1\mid\vecf{X}\sim F$, with symmetric $F(-x)=1-F(x)$.
When $X_1=-F^{-1}(\alpha)$, for any $X_2$, $\Pr(H_0 \mid \vecf{X})=\alpha$ exactly; this is indicated by the right vertical line in the figure.
For smaller values of $X_1$ (to the left of the vertical line), $\Pr(H_0 \mid \vecf{X})>\alpha$, whereas for larger $X_1$ (to the right), $\Pr(H_0 \mid \vecf{X})<\alpha$, as indicated in the figure.
For any $\vecf{\theta}$ with $\phi(\vecf{\theta})=c_0$, i.e., with $\theta_1=0$ (along the left vertical line), the rejection probability $\Pr(X_1>-F^{-1}(\alpha)\mid\vecf{\theta})=\alpha$ exactly, as indicated in the figure, so the test's size is also $\alpha$.
This is essentially the one-dimensional result for $X_1$ and $\theta_1$, applied to all $X_2$ and $\theta_2$.

\Cref{fig:thm-intuition-sub} illustrates \cref{thm:sub}.
The setup is the same as for \cref{fig:thm-intuition-half} but with a different null hypothesis $H_0 \colon \vecf{\theta}\in\Theta_0$ (the shaded region on the left).
The set $\Delta$ from \cref{thm:sub} is the white area between the shaded $\Theta_0$ and the vertical (dashed) line $\phi(\vecf{\theta})=c_0$, which again can be thought of as $\theta_1=0$.
Consider any point on the right vertical line, where previously $\Pr(H_0 \mid \vecf{X})=\alpha$ in \cref{fig:thm-intuition-half}.
Since $\Delta$ has positive posterior probability given any $\vecf{X}$ on that line, $H_0$ must have strictly smaller posterior probability than in \cref{fig:thm-intuition-half}, i.e., strictly smaller than $\alpha$.
Thus, the boundary of the rejection region must shift to the left by some amount, enlarging the rejection region.
Consequently, for any point on the left vertical line where previously the rejection probability was exactly $\alpha$, the rejection probability is now strictly above $\alpha$.
Most such points are now outside $\Theta_0$ (so rejection probability is power, not type I error rate), but there is one such point on the boundary of $\Theta_0$.

\Cref{fig:thm-intuition-sub} could be drawn for any support point of $\Theta_0$, by redefining $\phi(\cdot)$.
That is, consider any point on the boundary of $\Theta_0$ that is also on a supporting hyperplane (here, tangent line) $\phi(\vecf{\theta})=c_0$, meaning $\Theta_0$ is entirely on one side of the line.
The picture will look qualitatively similar at any such point.
The key is that after $\Delta$ is carved away from the supporting half-space $\phi(\vecf{\theta})\le c_0$, then all points with $\phi(\vecf{\theta})=c_0$ have rejection probability above $\alpha$, including the point on the boundary of $\Theta_0$.

\Cref{fig:thm-intuition-super} illustrates \cref{thm:super}.
The setup is again the same except the null hypothesis.
Now $\Theta_0$ is larger than the original half-space, so $\Delta$ is the shaded area within $\Theta_0$ to the right of the left vertical line ($\phi(\vecf{\theta})=c_0$).
Reversing the logic from before, this \emph{increases} the posterior probability of $H_0$ for all points $\vecf{X}$ on the right vertical line.
This in turn shifts the boundary of the rejection region to the right, making it smaller.
Thus, the rejection probability for any $\vecf{\theta}$ on the left vertical line is smaller than in \cref{fig:thm-intuition-half}, i.e., strictly below $\alpha$, including the point on the boundary of $\Theta_0$.
Again, by redefining $\phi(\cdot)$, essentially the same picture could be drawn for other points on the boundary of $\Theta_0$, specifically any support points.

\subsection{Discussion of results}
\label{sec:results-discussion}

\Cref{thm:sub} is partly a result of the prior $\Pr(H_0)$ being ``small'' when $\Theta_0$ is small.
That is, since the prior over the \emph{parameter} is always the same (regardless of $\Theta_0$), the implicit prior $\Pr(H_0)$ shrinks when $\Theta_0$ shrinks.
(Technically, the prior in \Cref{a:multi} is improper, so $\Pr(H_0)$ is not well-defined, but the qualitative idea remains.)
Unless $\Theta_0$ is a half-space, this differs from \citet{BergerSellke1987} and others who only consider ``objective'' priors with $\Pr(H_0)=1/2$.
Whether placing a prior on the null (like $\Pr(H_0)=1/2$) or on the parameter is more appropriate depends on the empirical setting.
Often it is easier computationally not to set a specific $\Pr(H_0)$; e.g., one may use the same posterior to compute probabilities of many different hypotheses.
However, this can lead to a small $\Pr(H_0)$ and thus large type I error rate.

Although the implicitly small $\Pr(H_0)$ partially explains \cref{thm:sub}, the shape of $\Theta_0$ also plays an important role.
For example, in $\R^2$, let $\vecf{\theta}=(\theta_1,\theta_2)$ and $H_0 \colon \theta_1\theta_2 \ge 0$, so $\Theta_0$ comprises the first and third quadrants (and thus is not contained in any half-space).
This $\Theta_0$ is the same ``size'' as the half-space $\{\vecf{\theta} : \theta_1 \ge 0\}$.
However, with bivariate normal sampling and posterior distributions, \cref{thm:half} implies the Bayesian test of the half-space has size $\alpha$, whereas the size of the Bayesian test of $H_0 \colon \theta_1\theta_2\ge0$ may be either strictly above or equal to $\alpha$, depending on the correlation.
For example, let $\vecf{X}=(X_1,X_2)$ have a bivariate normal sampling distribution with $\Corr(X_1,X_2)=-1$.
Then the test is equivalent to a scalar test where $H_0$ is a finite, closed interval, in which case the Bayesian test's size strictly exceeds $\alpha$ by \cref{thm:sub}.
The same holds for other strong negative correlations, although size decreases to $\alpha$ as the correlation increases; see \cref{sec:app-bivariate} as well as \citet[\S3.2]{Kaplan2015b} for details.
Also, simulation results in \cref{sec:ex-ordinal} show size can exceed $\alpha$ even with prior $\Pr(H_0)=1/2$.

As a special case, \cref{thm:sub} can apply to $\Theta_0=\{\vecf{\theta} : g(\vecf{\theta}) \le g_0\}$ when $g(\cdot)$ is directionally differentiable, as in \citet{FangSantos2018} and others.
For example, all the examples in Section 4 of \citet{FangSantos2018} concern $\vecf{\theta}$ belonging to a convex set.
They provide a frequentist resampling scheme that is consistent and corresponding hypothesis tests that control asymptotic size.
Thus, in cases like this where \cref{thm:sub} applies, not only is the Bayesian test's size strictly above $\alpha$, but it is also strictly above the size of an available frequentist test.

Many economic examples satisfy the condition of $\Theta_0$ being a subset of a half-space as in \cref{thm:subeq,thm:sub}.
Switching the null and alternative hypotheses then generates examples where $\Theta_0$ contains a half-space as in \cref{thm:non,thm:super}.
Often both are relevant; e.g., both the null of first-order stochastic dominance \citep[e.g.,][]{BarrettDonald2003} and the null of non-dominance \citep{DavidsonDuclos2013} have been studied.
For example, the null hypothesis of satisfaction of curvature constraints on economic functions (cost, production, etc.)\ has $\Theta_0$ smaller than a half-space, and thus the null of non-satisfaction has $\Theta_0$ larger than a half-space.
The ordinal distribution relationships in \cref{sec:ex-ordinal} and \citet{KaplanZhuo2019b} also have $\Theta_0$ smaller than a half-space.
The examples in Section 4 of \citet{FangSantos2018} also satisfy this condition, ``encompass[ing] tests of moment inequalities, shape restrictions, and the validity of random utility models'' (\S4, p.\ 25), the latter referring to \citet{KitamuraStoye2018}.
The general moment inequality null hypothesis $H_0 \colon \E[\vecf{W}] \le \vecf{0}$ with $\vecf{W} \in \R^d$ as in equation (58) of \citet{FangSantos2018} covers many other applications, but it has been studied already by \citet{Kline2011} (see \cref{sec:lit-close}); and although technically not covered by his results, first-order stochastic dominance with continuous distributions is essentially an infinite-dimensional extension of this setup.
Example 4.2 of \citet{FangSantos2018} considers shape restrictions on the infinite-dimensional regression quantile process, e.g., the restriction that the coefficient on the regressor of interest is monotonic in the quantile index.
%
%
%
In dynamic macroeconomic models, certain parameter values lead to a stable equilibrium while others do not.
Studies of such bifurcation find that $95\%$ confidence sets usually include values with an unstable equilibrium; more directly, equilibrium (in)stability could be tested.
Even in relatively simple models, the inequalities characterizing stability involve nonlinear functions of multiple parameters.
See for example \citet{BarnettDuzhak2010}, \citet{BarnettEryilmaz2016}, and references therein.
Finally, as described in the introduction, many general hypotheses about identified sets correspond to the point-identified parameter lying in some region $\Theta_0$.

Whether a set $\Theta$ is treated as the null or alternative hypothesis affects the properties of the Bayesian test.
Combining \cref{thm:sub,thm:non}, the Bayesian test of $H_0 \colon \vecf{\theta}\in\Theta$ may have size strictly above $\alpha$ while the Bayesian test of $H_0 \colon \vecf{\theta}\not\in\Theta$ has size strictly below or equal to $\alpha$.
Moreover, if $\vecf{\theta}$ is a support point of $\Theta$, then the rejection probability of $H_0\colon\vecf{\theta}\in\Theta$ is strictly above $\alpha$ by \cref{thm:sub}; but the rejection probability of $H_0\colon\vecf{\theta}\not\in\Theta$ is strictly below $\alpha$ by \cref{thm:super}.

Many nonlinear inequalities could be recast as linear inequalities, but at the expense of additional approximation error. 
For example, for $g \colon \R^d \mapsto \R$ and parameter $\vecf{\mu} \in \R^d$, the nonlinear $H_0 \colon g(\vecf{\mu}) \le 0$ could be written as $H_0 \colon \beta \le 0$ with $\beta \equiv g(\vecf{\mu})$. 
By the (first-order) delta method \citep[e.g.,][Thm.\ 6.8]{Hansen2020econometrics}, if $g(\cdot)$ is continuously differentiable in a neighborhood of $\vecf{\mu}$ and $\vecf{G} \equiv \left. \pD{}{\vecf{u}} g(\vecf{u}) \right\rvert_{\vecf{u}=\vecf{\mu}}$, then $\sqrt{n}(\hat{\vecf{\mu}}-\vecf{\mu}) \dconv \NormDist(\vecf{0}, \matf{V})$ implies $\sqrt{n}(\hat\beta-\beta) \dconv \NormDist(0, \vecf{G}'\matf{V}\vecf{G})$, where $\hat\beta \equiv g(\hat{\vecf{\mu}})$.
However, this approximation may be poor in finite samples.

To be concrete, consider the following delta method example, continuing the same notation.
Let $d=2$ and $\beta \equiv g(\vecf{\mu})=\mu_1^2+\mu_2^2-1$, so $\{\vecf{u} : g(\vecf{u}) \le 0\}$ is the unit disk in $\R^2$; the delta method gives $\sqrt{n}(g(\hat{\vecf{\mu}})-g(\vecf{\mu})) \dconv \NormDist(0,4\vecf{\mu}'\matf{V}\vecf{\mu})$.
Further imagine $\hat{\vecf{\mu}} \sim \NormDist(\vecf{\mu}, \matf{V}/n)$ in a finite sample of $n$ observations, with the corresponding posterior.
Then the Bayesian test's size is strictly above $\alpha$, as in \cref{thm:sub}.
Further, as the variance $\matf{V}/n$ grows, $\Pr(H_0\mid\hat{\vecf{\mu}})\to0$ given any $\hat{\vecf{\mu}}$, so size grows to $1$.
In apparent contradiction, $H_0 \colon \beta \le 0$ suggests the Bayesian test has approximately correct size, by \cref{thm:half}.
The ``contradiction'' is simply that the asymptotic approximation is less accurate due to the delta method's linear approximation of $g(\cdot)$.
We avoid this layer of delta method approximation by treating nonlinear inequalities directly, providing better finite-sample insights.
Another drawback of the delta method is the existence of singularity points; e.g., with $\vecf{\mu}=\vecf{0}$ in this example, the delta method's distribution of $g(\hat{\vecf{\mu}})$ is degenerate even if the distribution of $\hat{\vecf{\mu}}$ is not.

\subsection{Special cases of results}
\label{sec:results-special}

In the special case when $X,\theta \in \R$, similar results to \cref{thm:half} are found in the literature, like in \citet{CasellaBerger1987}. 
Less general versions of \cref{thm:subeq,thm:sub} have also been given when $X,\theta \in \R$.

A special case of \cref{thm:sub} with $\vecf{\theta}\in\R^d$ and $H_0 \colon \vecf{\theta}\ge\vecf{0}$ was explored by \citet{Kline2011}, as detailed in \cref{sec:lit-close}.

\Cref{thm:1} also includes the special case of linear inequalities in $\R^d$. 
\Cref{thm:half} states that for a single linear inequality $H_0 \colon \vecf{c}'\vecf{\theta} \le c_0$, the Bayesian test has size $\alpha$. 
\Cref{thm:sub} states that for multiple linear inequalities, the Bayesian test's size is strictly above $\alpha$, and its RP is strictly above $\alpha$ at every boundary point of $\Theta_0$.

\section{Examples}
\label{sec:ex}

We illustrate \cref{thm:1} through two specific examples; another is in \cref{sec:app-SD1}.

\subsection{Example: ordinal distributions}
\label{sec:ex-ordinal}

\subsubsection{Background}

The original motivation for this paper was inference on relationships between ordinal distributions.
\Citet{KaplanZhuo2019b} connect such ordinal relationships to relationships between corresponding latent distributions.
Here, we focus on one such result.

Consider the ordinal variable self-reported health status (SRHS), commonly used in health economics.
Let random variables $X$ and $Y$ denote SRHS from two different populations, like different socioeconomic or demographic groups.
Let $X=1$ denote ``poor'' health (the worst category), $X=5$ ``excellent'' health (the best), and $X=2,3,4$ the ordered intermediate categories (fair, good, very good), and similarly for $Y$.
These numbers have no cardinal meaning; they are only for notational convenience.
Imagine they are determined by an individual's underlying latent health $X^*$ or $Y^*$, where $X=1$ is reported if $X^*\le\gamma_1$, $X=2$ is reported if $\gamma_1<X^*\le\gamma_2$, etc.
The ordinal CDFs are $F_X(\cdot)$ and $F_Y(\cdot)$; the latent CDFs are $F_X^*(\cdot)$ and $F_Y^*(\cdot)$.

If $Y$ values are determined from $Y^*$ using the same $\gamma_j$ thresholds as $X$, possibly shifted by a common constant ($\gamma_j+c$), then it is possible to learn about the relative dispersion of latent $X^*$ and $Y^*$ from the ordinal $X$ and $Y$.
Specifically, if there is a ``single crossing'' of $F_Y(\cdot)$ by $F_X(\cdot)$ from below, then certain interquantile ranges of $X^*$ are smaller than those of $Y^*$.
(Other interquantile ranges lack evidence for either direction.)
``Single crossing'' (SC) means that $F_X(\cdot)$ is below $F_Y(\cdot)$ up to a certain category, and then is above thereafter.
That is, for some category $m$, $F_X(j)<F_Y(j)$ for $j<m$ and $F_X(j)>F_Y(j)$ for $j>m$.

This can also be interpreted in terms of partial identification.
The two ordinal CDFs are point-identified, but the two latent CDFs are only partially identified.
That is, there is a mapping from a pair of ordinal CDFs to a (very large) set of pairs of latent CDFs that are consistent with the ordinal pair.
The ordinal SC implies that \emph{all} pairs of latent CDFs in the identified set have certain interquantile ranges of $X^*$ smaller than those of $Y^*$.

\subsubsection{Connections to \texorpdfstring{\cref{thm:1}}{Theorem \ref{thm:1}}}

The null hypothesis of SC (of $Y$ by $X$) can be written in terms of the parameters as follows.
There are eight unknown ordinal parameters: $F_X(j)$ and $F_Y(j)$ for $j=1,2,3,4$, since $F_X(5)=F_Y(5)=1$.
Define the four CDF difference parameters
\begin{equation}\label{eqn:delta}
\theta_j \equiv F_X(j)-F_Y(j) .
\end{equation}
SC requires $\theta_1<0$ and $\theta_4>0$, but there are multiple possibilities for $\theta_2$ and $\theta_3$.
Specifically, SC rules out $\theta_2>0>\theta_3$, but the other three quadrants of $(\theta_2,\theta_3)$ satisfy SC.
We ignore the difference between $<$ and $\le$ since they are statistically indistinguishable.

Writing the SC null hypothesis as $H_0 \colon \vecf{\theta}\in\Theta_0$, the shape of $\Theta_0$ has both convex and non-convex features.
First, consider the first and fourth components.
Since $\theta_1<0<\theta_4$ is a necessary (not sufficient) condition,
\begin{equation}
\Theta_0 
\subset \{ \vecf{\theta} : \theta_1<0<\theta_4 \}
\subset \{ \vecf{\theta} : \theta_1 < 0 \} ,
\end{equation}
where $\phi(\vecf{\theta})=\theta_1$ is a continuous linear functional of $\vecf{\theta}$.
That is, $\Theta_0$ is a subset of a quarter-space, which is a subset of a half-space that satisfies the conditions in \cref{thm:sub}.
Thus, the Bayesian test's size should be strictly above $\alpha$, with rejection probability above $\alpha$ for DGPs with $\vecf{\theta}$ on the boundary of $\Theta_0$ near the corner $(\theta_1,\theta_4)=(0,0)$.

Second, looking at only $(\theta_2,\theta_3)$ reveals local non-convexity.
Imagine $\theta_1<0<\theta_4$ by a large margin, and the sample size is large enough that the posterior of $\theta_1<0<\theta_4$ is near one.
Then, the test practically only relies on $(\theta_2,\theta_3)$.
But $\Theta_0$ only excludes one quadrant of the $(\theta_2,\theta_3)$ subspace, the one with $\theta_2>0>\theta_3$.
When the conditions on $(\theta_1,\theta_4)$ are satisfied, this leaves three quadrants of the $(\theta_2,\theta_3)$ subspace, which is more than a half-space.
In this case, the Bayesian test's size probably still equals $\alpha$, but it may reject much less under certain DGPs, like with $(\theta_2,\theta_3)$ near the origin.
That is, \cref{thm:super} suggests asymptotic type I error rates strictly below $\alpha$ for DGPs local to $(\theta_2,\theta_3)=(0,0)$ with $\theta_1<0<\theta_4$, but since $\Theta_0$ is not larger than a half-space globally, this may not be true in small samples.

Instead of SC, the null hypothesis could be the opposite, i.e., non-SC.
Then, when $(\theta_1,\theta_4)=(0,0)$, instead of \cref{thm:sub} implying type I error rate above $\alpha$, \cref{thm:super} implies type I error rate below $\alpha$.
When $(\theta_2,\theta_3)=(0,0)$ and $\theta_1<0<\theta_4$, instead of \cref{thm:super} applying asymptotically, \cref{thm:sub} suggests asymptotic type I error rate above $\alpha$.

\subsubsection{Simulation setup}

The following simulations show the importance of the shape of $\Theta_0$ and sometimes sample size.
Independent, iid samples of $X$ and $Y$ are taken, with sample size $n$ for both.
The results show the simulated rejection probability (RP) of the Bayesian test of $X$ single crossing $Y$.
DGP 1 is at the convex ``corner'' with $(\theta_1,\theta_4)=(0,0)$, where the Bayesian test's RP should be above $\alpha$.
Specifically, $F_X(1)=F_Y(1)=0.2$, $F_X(2)=0.39<0.4=F_Y(2)$, $F_X(3)=0.5<0.6=F_Y(3)$, and $F_X(4)=F_Y(4)=0.8$.
DGP 2 is at the local non-convexity with $(\theta_2,\theta_3)=(0,0)$, where the RP may be below $\alpha$ with large enough $n$.
Specifically, $F_X(1)=0.18<0.22=F_Y(1)$, $F_X(2)=F_Y(2)=0.4$, $F_X(3)=F_Y(3)=0.6$, and $F_X(4)=0.82>0.78=F_Y(4)$.
Since both DGPs are on the boundary of SC and non-SC, the null hypothesis can be considered either SC or non-SC, and the rejection probability can be interpreted as a type I error rate.
Each simulation had $1000$ replications, with $1000$ posterior draws each.

The Bayesian test is computed as follows.
Let $\vecf{p}_X$ be the vector with entry $j$ equal to $\Pr(X=j)$, $j=1,\ldots,5$, and similarly for $\vecf{p}_Y$.
Let $\vecf{n}_X$ be the vector of data counts with entry $j$ equal to the number of observations with $X_i=j$ (out of $i=1,\ldots,n$), and similarly for $\vecf{n}_Y$.
With a uniform prior, the posterior distributions are Dirichlet; specifically, with $\vecf{1}\equiv(1,1,1,1,1)$,
\begin{equation}
\vecf{p}_X \mid \vecf{n}_X \sim \DirDist(\vecf{n}_X+\vecf{1}) , \quad
\vecf{p}_Y \mid \vecf{n}_Y \sim \DirDist(\vecf{n}_Y+\vecf{1}) .
\end{equation}
The Bayesian test takes draws from these posteriors and computes the proportion of draws in which SC holds, rejecting the SC null hypothesis if the proportion is below $\alpha$, and rejecting the null of non-SC if the proportion is above $1-\alpha$.

Additionally, to help isolate the role of ``shape'' from ``size'' of $\Theta_0$, a test with adjusted prior is computed.
It is identical except that the prior is adjusted to have $1/2$ probability on $H_0$.
Here, the original probability of $H_0$ can be simulated by drawing from the uniform prior, and then the prior is normalized using this original probability.
The prior remains uniform within $\Theta_0$, and uniform within $\Theta_0^\complement$, but at two different levels instead of the same level.
See \citet{GoutisEtAl1996} or \citet{KaplanZhuo2019b} for details.

\begin{table}[htbp]
\centering
\caption{\label{tab:sim-ordinal}Simulated type I error rate of Bayesian tests of ordinal single crossing, $\alpha=0.1$.}
\begin{tabular}[c]{S[table-format=1.0,round-precision=0]S[table-format=3.0,round-mode=none] S[table-format=1.3,round-precision=3,round-mode=places]S[table-format=1.3,round-precision=3,round-mode=places] c S[table-format=1.3,round-precision=3,round-mode=places]S[table-format=1.3,round-precision=3,round-mode=places]}
\toprule
&&\multicolumn{2}{c}{$H_0 \colon $SC}
&&\multicolumn{2}{c}{$H_0 \colon$non-SC} \\
\cmidrule{3-4}\cmidrule{6-7}
{DGP} & {$n$} & {RP} & {RP (adj)} && {RP} & {RP (adj)} \\
\midrule
1 &   20 & 0.422 & 0.170 && 0.000 & 0.014 \\
1 &  100 & 0.484 & 0.186 && 0.000 & 0.010 \\
1 &  500 & 0.531 & 0.249 && 0.000 & 0.008 \\[2pt]
2 &   20 & 0.246 & 0.056 && 0.001 & 0.050 \\
2 &  100 & 0.128 & 0.022 && 0.009 & 0.148 \\
2 &  500 & 0.011 & 0.002 && 0.130 & 0.501 \\
\bottomrule
\end{tabular}
\end{table}

\subsubsection{Simulation results}

\Cref{tab:sim-ordinal} shows the results.
``RP'' is the rejection probability of the first Bayesian test, and ``RP (adj)'' is the rejection probability of the Bayesian test with the adjusted prior.
Some patterns emerge.

For DGP 1, RP is indeed well above $\alpha$ when the null hypothesis is SC, consistent with \cref{thm:sub}.
Here, the RP actually increases as $n$ increases; it is not simply a finite-sample difference.
Even with the adjusted prior, RP is well above $\alpha$ and again increasing with $n$.
In contrast, when the null is non-SC, the RP is near zero, consistent with \cref{thm:super}.

For DGP 2, the sample size plays an important role.
When the sample size is small, the posterior and sampling distributions are very dispersed, so the global shape of $\Theta_0$ matters most.
Globally, for the null of SC, $\Theta_0$ is a subset of a quarter-space; \cref{thm:sub} suggests DGPs near the boundary may have RP above $\alpha$.
Indeed, RP is again well above $\alpha$ for small $n$.
With the null of non-SC, instead \cref{thm:super} suggests RP below $\alpha$, which is seen in the table.
As $n$ increases, the local shape of $\Theta_0$ matters more, so \cref{thm:super} is more appropriate for the null of SC and \cref{thm:sub} is more appropriate for the null of non-SC.
In this case, with the largest sample size, RP is far below $\alpha$ for the SC null, opposite the result with small $n$.
Similarly, the RP becomes above $\alpha$ for the non-SC null, also opposite the small $n$ results.
Like with DGP 1, the RP with adjusted prior is different, but it is not always closer to $\alpha$.

\subsection{Example: curvature constraints}
\label{sec:ex-curvature}

\subsubsection{Setting and connections to \texorpdfstring{\cref{thm:1}}{Theorem \ref{thm:1}}}

One common nonlinear inequality hypothesis in economics is a curvature constraint like concavity.
Such constraints come from economic theory, often the second-order condition of an optimization problem like utility maximization or cost minimization.
As noted by \citet{ODonnellCoelli2005}, the Bayesian approach is appealing for imposing or testing curvature constraints due to its (relative) simplicity.
However, according to \cref{thm:1}, since curvature is usually satisfied in a parameter subspace much smaller than a half-space, Bayesian inference similar to \cref{meth:Bayes} may be much less favorable toward the curvature hypothesis than frequentist inference would be; i.e., the size of the Bayesian test in \cref{meth:Bayes} may be well above $\alpha$. 
This is true in \cref{tab:sim-translog} below.

Our example concerns concavity of a cost function with the ``translog'' functional form \citep{ChristensenEtAl1973}. 
This has been a popular way to parameterize cost, 
indirect utility, 
and production functions. 
The translog is more flexible than many traditional functional forms, allowing violation of certain implications of economic theory (such as curvature) without reducing such constraints to the value of a single parameter. 
Since \citet{Lau1978}, there has been continued interest in methods to impose curvature constraints during estimation, as well as methods to test such constraints. 
Although ``flexible,'' the translog is still parametric, so violation of curvature constraints may come from misspecification of the functional form rather than violation of economic theory.%
\footnote{With a nonparametric model, one may more plausibly test the theory itself, although there are always other assumptions that may be violated; see \citet{DetteEtAl2016} for nonparametrically testing negative semidefiniteness of the Slutsky substitution matrix.}

Specifically, we test concavity of cost in input prices as follows.\footnote{The ``translog'' example on page 346 of \citet{Dufour1989} is even simpler but appears to ignore the fact that second derivatives are not invariant to log transformations.} 
With output $y$, input prices $\vecf{w}=(w_1,w_2,w_3)$, and total cost $C(y,\vecf{w})$, the translog model is 
\begin{equation}\label{eqn:translog}
\begin{split}
\ln(C(y,\vecf{w})) 
= a_0 &+ a_y\ln(y) + (1/2) a_{yy} [\ln(y)]^2 
 +\sum_{k=1}^{3}a_{yk}\ln(y)\ln(w_k)
\\&
 +\sum_{k=1}^{3} b_k\ln(w_k) 
 +(1/2)\sum_{k=1}^{3}\sum_{m=1}^{3} b_{km}\ln(w_k)\ln(w_m) . 
\end{split}
\end{equation}
Standard economic assumptions imply that $C(y,\vecf{w})$ is concave in $\vecf{w}$ \citep[as in][\S7.3]{Kreps1990}, which corresponds to the Hessian matrix (of $C$ with respect to $\vecf{w}$) being negative semidefinite (NSD), which in turn corresponds to all the Hessian's principal minors of order $p$ (for all $p=1,2,3$) having the same sign as $(-1)^p$ or zero.\footnote{In some cases, only leading principal minors are required; see \citet{Mandy2018}.}

For simplicity, we consider local concavity at the point $(1,1,1,1)$: 
\begin{equation}\label{eqn:H0-translog-local}
H_0 \colon  \matf{H} \equiv \left.\frac{\partial^2 C(y,\vecf{w})}{\partial\vecf{w}\partial\vecf{w}'}\right|_{(y,\vecf{w})=(1,1,1,1)} 
\textrm{ is NSD} . 
\end{equation}
This is necessary but not sufficient for global concavity; rejecting local concavity implies rejection of global concavity. 
In \cref{sec:app-translog}, we show that even this weaker constraint corresponds to a set of parameter values much smaller than a half-space, so \cref{thm:sub} applies.

\subsubsection{Simulation setup and results}

Our simulation DGP is as follows. 
To impose homogeneity of degree one in input prices, we use the normalized model (with error term $\epsilon$ added) 
\begin{equation}\label{eqn:translog-normalized}
\begin{split}
\ln(C/w_3) &= a_0 + a_y\ln(y) +(1/2)a_{yy}[\ln(y)]^2
            +\sum_{k=1}^{2}a_{yk}\ln(y)\ln(w_k/w_3)
\\&\quad
            +\sum_{k=1}^{2}b_k\ln(w_k/w_3)
            +(1/2)\sum_{k=1}^{2}\sum_{m=1}^{2}b_{km}\ln(w_k/w_3)\ln(w_m/w_3)
            +\epsilon 
\end{split}
\end{equation}
for both data generation and inference.%
\footnote{Alternatively, cost share equations may be used. 
Shephard's lemma implies that the demand for input $k$ is $x_k=\partial C/\partial w_k$. 
The cost share for input $k$ is then 
$s_k = x_k w_k / C
     = (\partial C/\partial w_k) (w_k/C)
     = \partial \ln(C)/\partial \ln(w_k) 
= b_k + a_{yk}\ln(y) + \sum_{j=1}^{3}b_{jk}\ln(w_j) $. %
} 
The parameter values are $b_1=b_2=1/3$, $b_{11}=b_{22}=2/9-\delta$ (more on $\delta$ below), and $b_{12}=-1/9$ to make some of the inequality constraints in $H_0$ close to binding, as well as $a_0=1$, $a_y=1$, $a_{yy}=0$, $a_{yk}=0$. 
The other parameter values follow from imposing symmetry ($b_{km}=b_{mk}$) and homogeneity. 
When $\delta=0$, $\matf{H}$ is a matrix of zeros, on the boundary of being NSD in that each principal minor equals zero. 
When $\delta>0$, all principal minors satisfy NSD and are non-zero (other than $\det(\matf{H})=0$, which is always true under homogeneity). 
We set $\delta=0.001$. 
In each simulation replication, an iid sample is drawn, where $\ln(y)$ and all $\ln(w_k)$ are $\NormDistp{0}{\sigma=0.1}$, $\epsilon\sim\NormDistp{0}{\sigma_\epsilon}$, and all variables are mutually independent. 
There are $n=\num{100}$ observations per sample, $\num{500}$ simulation replications, and $\num{200}$ posterior draws per replication. 
The local monotonicity constraints $b_1,b_2,b_3\ge0$ were satisfied in $100.0\%$ of draws among all replications.

\Cref{tab:sim-translog} reports values from two methods. 
For the method denoted ``Bayesian bootstrap'' in the table header, the posterior probability of $H_0$ is computed by a nonparametric Bayesian method with improper Dirichlet process prior, i.e., the Bayesian bootstrap of \citet{Rubin1981} based on \citet{Ferguson1973} and more recently advocated in economics by \citet{ChamberlainImbens2003}. 
For the method denoted ``Normal,'' the parameter vector is sampled from a normal distribution with mean equal to the ordinary least squares (OLS) estimate and covariance matrix equal to the corresponding (homoskedastic) asymptotic covariance matrix estimate; this is the posterior from a homoskedastic normal linear regression model with improper prior (or asymptotically). 
To accommodate numerical imprecision, we deem an inequality satisfied if it is within $10^{-7}$. 
The simulated type I error rate is the proportion of simulated samples for which the posterior probability of $H_0$ was below $\alpha$.

\begin{table}[htbp]
\centering
\caption{\label{tab:sim-translog}Simulated type I error rate of Bayesian tests of local NSD.}
\begin{tabular}[c]{S[table-format=1.2,round-precision=2,round-mode=places]S[table-format=1.2,round-precision=2,round-mode=places] S[table-format=1.3,round-precision=3,round-mode=places]S[table-format=1.3,round-precision=3,round-mode=places]}
\toprule
           &                     & {Bayesian}  &  \\
{$\alpha$} & {$\sigma_\epsilon$} & {bootstrap} & {Normal} \\
\midrule
0.05 & 0.00 & 0.000 & 0.000 \\
0.05 & 0.10 & 0.112 & 0.084 \\
0.05 & 0.20 & 0.354 & 0.318 \\
0.05 & 0.30 & 0.554 & 0.532 \\
0.05 & 0.40 & 0.676 & 0.694 \\
0.05 & 0.50 & 0.764 & 0.772 \\[2pt]
0.10 & 0.00 & 0.000 & 0.000 \\
0.10 & 0.10 & 0.186 & 0.160 \\
0.10 & 0.20 & 0.530 & 0.526 \\
0.10 & 0.30 & 0.740 & 0.764 \\
0.10 & 0.40 & 0.844 & 0.872 \\
0.10 & 0.50 & 0.890 & 0.910 \\
\bottomrule
\end{tabular}
\end{table}

\Cref{tab:sim-translog} shows the type I error rate of the Bayesian tests of \cref{eqn:H0-translog-local} given our DGP. 
The values of $\alpha$ and $\sigma_\epsilon$ are varied as shown in the table. 
The two Bayesian tests are very similar, always within a few percentage points of each other. 
As a sanity check, when $\sigma_\epsilon=0$, the RP is zero since the constraints are satisfied by construction. 
As $\sigma_\epsilon$ increases, the RP increases well above $\alpha$, even over $50\%$.%
\footnote{The results with $\delta=0.01$ are similar; with $\delta=0$, RP jumps to over $80\%$ even with $\sigma_\epsilon=0.001$.} 
This is consistent with \cref{thm:sub}.


\section{Conclusion}
\label{sec:conclusion}

We have explored the frequentist properties of Bayesian inference on general nonlinear inequality constraints, including hypotheses about identified set properties, providing formal results on the role of the shape of the null hypothesis parameter subspace.
Future work could include investigation of approaches like \citet{MuellerNorets2016} applied to nonlinear inequality testing.
Moreover, one could explore how to match frequentist size by adjusting the prior or loss function.

\singlespacing
\bibliographystyle{jpe}

\onehalfspacing

\appendix

\numberwithin{equation}{section}

\section{Proofs}
\label{sec:app-pfs}

\begin{proof}[\bfseries Proof of \cref{thm:1}]
For \cref{thm:half}: 
the Bayesian test rejects iff 
\begin{equation*}
\alpha 
\ge \Pr\bigl(\phi(\vecf{\theta})\le c_0\mid \vecf{X} \bigr)
= \Pr\bigl( \phi(\vecf{\theta})-\phi(\vecf{X}) \le c_0-\phi(\vecf{X}) \mid \vecf{X} \bigr)
\equiv F\bigl(c_0-\phi(\vecf{X})\bigr) . 
\end{equation*}
Given any $\vecf{\theta}$ such that $\phi(\vecf{\theta})\le c_0$ (so $H_0$ holds), the RP is 
\begin{align*}
\Pr\bigl( F(c_0-\phi(\vecf{X})) \le \alpha \mid \vecf{\theta} \bigr) 
&= \Pr\bigl( \overbrace{1-F(\phi(\vecf{X})-c_0)}^{\textrm{by symmetry}} \le \alpha \mid \vecf{\theta} \bigr) 
\\&= \Pr\bigl( F(\phi(\vecf{X})-c_0) \ge 1-\alpha \mid \vecf{\theta} \bigr) 
\\&\overbrace{\le \Pr\bigl( F(\phi(\vecf{X})-\phi(\vecf{\theta})) \ge 1-\alpha \mid \vecf{\theta} \bigr) }^{\textrm{since $\phi(\vecf{\theta})\le c_0$ under $H_0$}}
\\&= \alpha 
\end{align*}
since $F\bigl(\phi(\vecf{X})-\phi(\vecf{\theta})\bigr) \mid \vecf{\theta} \sim \UnifDist(0,1)$. 
If $\phi(\vecf{\theta})=c_0$, then the $\le$ becomes $=$.

For \cref{thm:subeq}: 
because $\Theta_0\subseteq\{\vecf{\theta}:\phi(\vecf{\theta})\le c_0\}$, then for any $\vecf{X}$, 
\begin{equation*}
\Pr(\vecf{\theta}\in\Theta_0\mid \vecf{X}) \le \Pr\bigl( \phi(\vecf{\theta})\le c_0\mid \vecf{X} \bigr) . 
\end{equation*}
Consequently, the rejection region for $H_0 \colon \vecf{\theta}\in\Theta_0$ is at least as big as the rejection region for $H_0 \colon \phi(\vecf{\theta})\le c_0$: for some $r\in\R$, 
\begin{equation}\label{eqn:subeq-region}
\begin{split}
& \mathcal{R}_1 \subseteq \mathcal{R}_2, \;
\mathcal{R}_1 \equiv \{ \vecf{X} : \Pr(\phi(\vecf{\theta})\le c_0 \mid \vecf{X}) \le \alpha \} 
= \{\vecf{X} : \phi(\vecf{X})\ge r\}, \\
& \mathcal{R}_2 \equiv \{ \vecf{X} : \Pr(\vecf{\theta}\in\Theta_0\mid \vecf{X}) \le \alpha \} . 
\end{split}
\end{equation}
Given any $\vecf{\theta}\in\Theta_0$, the probability that $\vecf{X}$ falls in the new, larger rejection region ($\mathcal{R}_2$) is at least as big as the probability that $\vecf{X}$ falls in the old, smaller rejection region ($\mathcal{R}_1$) from \cref{thm:half}. 
In particular, when $\phi(\vecf{\theta})=c_0$, the RP was exactly $\alpha$ in \cref{thm:half}. 
Since the new rejection region is weakly larger, the new RP when $\phi(\vecf{\theta})=c_0$ must be at least $\alpha$. 
If $c_0\in\phi(\Theta_0)$, then the proof is complete. 
Otherwise, with $\mathcal{R}_1$ and $\mathcal{R}_2$ from \cref{eqn:subeq-region}, and $\vecf{\theta}^*$ any value such that $\phi(\vecf{\theta}^*)=c_0$ (with the limit formed by a sequence of $\vecf{\theta}$ within $\Theta_0$), 
\begin{align}\notag
& \sup_{\vecf{\theta}\in\Theta_0} \Pr(\vecf{X}\in\mathcal{R}_2\mid\vecf{\theta}) \overbrace{\ge \lim_{\vecf{\theta}\to \vecf{\theta}^*} \Pr(\vecf{X}\in\mathcal{R}_2\mid\vecf{\theta})}^{\textrm{since $c_0\in\phi(\closure{\Theta_0})$}}
\overbrace{\ge \lim_{\vecf{\theta}\to \vecf{\theta}^*} \Pr(\vecf{X}\in\mathcal{R}_1\mid\vecf{\theta})}^{\textrm{by \cref{eqn:subeq-region}}}
\\&= \lim_{\vecf{\theta}\to \vecf{\theta}^*} \Pr\bigl( \phi(\vecf{X})\ge r \mid \vecf{\theta} \bigr) 
= \lim_{\vecf{\theta}\to \vecf{\theta}^*} \Pr(\phi(\vecf{X})-\phi(\vecf{\theta}) \ge r-\phi(\vecf{\theta}) \mid \vecf{\theta})
= \lim_{\vecf{\theta}\to \vecf{\theta}^*} 1-F\bigl( r - \phi(\vecf{\theta}) \bigr)
\notag 
\\&\overbrace{= 1-F(r-c_0)}^{\textrm{by continuity of $F,\phi$}} 
\overbrace{= \alpha}^{\textrm{by \cref{thm:half}}} . 
\label{eqn:subeq-size}
\end{align}

For \cref{thm:sub}: 
given the stated assumption that the posterior distribution of $\phi_2(\vecf{\theta})$ has a strictly positive PDF for any $\phi_2(\vecf{X})$, and the assumption that $\Delta$ has positive Lebesgue measure, then 
\begin{equation}\label{eqn:Delta-prob}
\Pr(\phi_2(\vecf{\theta}) \in \Delta \mid \phi_2(\vecf{X})) > 0\textrm{ for any }\phi_2(\vecf{X}). 
\end{equation}
Similar to \cref{eqn:subeq-region}, for some $r \in \R$, let 
\begin{align}\notag
& \mathcal{R}_a \subseteq \mathcal{R}_b, \;
\mathcal{R}_a \equiv \{ \phi_2(\vecf{X}) : \Pr(\phi_3(\phi_2(\vecf{\theta})) \le c_0 \mid \phi_2(\vecf{X})) \le \alpha \} 
= \{\phi_2(\vecf{X}) : \phi_3(\phi_2(\vecf{X})) \ge r\}, \\
& \mathcal{R}_b \equiv \{ \phi_2(\vecf{X}) : \Pr(\phi_2(\vecf{\theta}) \in \Phi_2 \mid \phi_2(\vecf{X})) \le \alpha \} . 
\label{eqn:sub-region}
\end{align}
Let $\vecf{p}^*$ be any value such that $\phi_3(\vecf{p}^*)=r$. 
Then, 
\begin{equation*}
\begin{split}
\Pr(\vecf{\theta} \in \Theta_0 \mid \phi_2(\vecf{X})=\vecf{p}^*) 
& \overbrace{\le \Pr( \phi_2(\vecf{\theta}) \in \Phi_2 \mid \phi_2(\vecf{X})=\vecf{p}^*)}^{\textrm{since }\vecf{\theta} \in \Theta_0  \implies  \phi_2(\vecf{\theta}) \in \Phi_2}
\\&= \overbrace{\Pr(\phi_3(\phi_2(\vecf{\theta})) \le c_0 \mid \phi_2(\vecf{X})=\vecf{p}^*)}^{=\alpha\textrm{ by \cref{thm:half}}}
-\overbrace{\Pr(\phi_2(\vecf{\theta}) \in \Delta \mid \phi_2(\vecf{X})=\vecf{p}^*)}^{>0\textrm{ by \cref{eqn:Delta-prob}}} 
\\& < \alpha . 
\end{split}
\end{equation*}
By the assumption that $\Pr(\phi_2(\vecf{\theta}) \in \Phi_2 \mid \phi_2(\vecf{X}))$ is continuous in $\phi_2(\vecf{X})$, there is some $\epsilon$-ball $\mathcal{B}$ around $\vecf{p}^*$ for which $\Pr(\vecf{\theta} \in \Theta_0 \mid \phi_2(\vecf{X}) \in \mathcal{B}) < \alpha$, too. 
The ball $\mathcal{B}$ has positive Lebesgue measure, as does the part of it lying outside $\mathcal{R}_a$ (i.e., $\mathcal{B}\cap\mathcal{R}_a^\complement$) since $\vecf{p}^*$ is on the boundary of $\mathcal{R}_a$. 
Since the sampling distribution of $\phi_2(\vecf{X})$ given any $\vecf{\theta}$ has a strictly positive PDF (by assumption),
\begin{equation}\label{eqn:ball-prob}
\Pr\bigl( \phi_2(\vecf{X}) \in (\mathcal{B}\cap\mathcal{R}_a^\complement) \mid \vecf{\theta} \bigr) > 0\textrm{ for any }\vecf{\theta} . 
\end{equation}
Also, using the assumption in \cref{a:multi} that the distribution of $\phi_2(\vecf{\theta})$ only depends on $\vecf{X}$ through $\phi_2(\vecf{X})$, as well as the assumption that $\vecf{\theta}\in\Theta_0 \implies \phi_2(\vecf{\theta})\in\Phi_2$, 
\begin{equation*}
    \Pr( \phi_2(\vecf{\theta}) \in \Phi_2 \mid \phi_2(\vecf{X}) )
  = \Pr( \phi_2(\vecf{\theta}) \in \Phi_2 \mid \vecf{X} )
\ge \Pr( \vecf{\theta} \in \Theta_0 \mid \vecf{X} ) , 
\end{equation*}
so $\Pr( \phi_2(\vecf{\theta}) \in \Phi_2 \mid \phi_2(\vecf{X}) ) \le \alpha  \implies  \Pr( \vecf{\theta} \in \Theta_0 \mid \vecf{X} ) \le \alpha$. 
Consequently, 
\begin{equation}\label{eqn:Rb-implies-R2}
\phi_2(\vecf{X}) \in \mathcal{R}_b  \implies  \vecf{X} \in \mathcal{R}_2 . 
\end{equation}
Letting $\vecf{\theta}^*$ be a value such that $\phi_3(\phi_2(\vecf{\theta}^*))=\phi(\vecf{\theta}^*)=c_0$, the Bayesian test's rejection probability is bounded from below by 
\begin{align*}
\Pr(\vecf{X} \in \mathcal{R}_2 \mid \vecf{\theta}^*)
&\overbrace{%
\ge \Pr(\phi_2(\vecf{X}) \in \mathcal{R}_b \mid \vecf{\theta}^*)
}^{\textrm{by \cref{eqn:Rb-implies-R2}}}
\\&\ge \overbrace{\Pr(\phi_2(\vecf{X}) \in \mathcal{R}_a \mid \vecf{\theta}^*)}^{=\alpha\textrm{ by \cref{thm:half}}}
    +\overbrace{\Pr\bigl( \phi_2(\vecf{X}) \in (\mathcal{B}\cap\mathcal{R}_a^\complement) \mid \vecf{\theta}^* \bigr) }^{>0\textrm{ by \cref{eqn:ball-prob}}}
\\&> \alpha . 
\end{align*}
As in the proof of \cref{thm:subeq}, 
if $\vecf{\theta}^* \in \Theta_0$, 
then the test's size is bounded below by $\Pr(\vecf{X} \in \mathcal{R}_2 \mid \vecf{\theta}^*)$ and the proof is complete. 
Otherwise, as before, the assumed continuity and $c_0 \in \phi(\closure{\Theta_0})$ imply 
\begin{align*}
\sup_{\vecf{\theta}\in\Theta_0} \Pr(\vecf{X} \in \mathcal{R}_2 \mid \vecf{\theta})
  &\ge 
\lim_{\vecf{\theta}\to\vecf{\theta}^*} \Pr( \vecf{X} \in \mathcal{R}_2 \mid \vecf{\theta} )
= \Pr( \vecf{X} \in \mathcal{R}_2 \mid \vecf{\theta}^* ) . 
\end{align*}

For \cref{thm:non}, some examples suffice. 
First, consider $H_0 \colon \phi(\vecf{\theta}) \ne 0$. 
Given any $\vecf{X}$, 
$\Pr(\phi(\vecf{\theta}) \ne 0  \mid \vecf{X})=1$ 
since $F$ is continuous, so the Bayesian test never rejects and its size is zero. 
Thus, size may be strictly below $\alpha$. 
Second, consider $H_0 \colon \phi(\vecf{\theta})\in\mathbb{Z}$ (the integers). 
This $H_0$ has zero posterior probability given any $\vecf{X}$, so the Bayesian test always rejects and has size equal to one. 
Thus, size may be strictly above $\alpha$. 
Third, consider the bivariate example with $H_0 \colon \theta_1 \le 0\textrm{ or }\theta_2 \le 0$. 
Let the sampling/posterior distribution be bivariate normal with unit variances and known correlation $\rho$. 
If $\rho=1$, then $X_1-\theta_1=X_2-\theta_2$ in every draw of $(X_1,X_2)$, i.e., it becomes a one-dimensional problem. 
Consequently, if $\theta_1-\theta_2 \ge 0$, then $X_1-X_2 \ge 0$ in every draw, and we may simply test $H_0 \colon \theta_2 \le 0$ using $X_2 \sim \NormDist(\theta_2,1)$, the Bayesian test of which has size $\alpha$. 
Similarly, if $\theta_1-\theta_2 \le 0$, then $X_1-X_2 \le 0$ in every draw, and the Bayesian test of $H_0 \colon \theta_1 \le 0$ has size $\alpha$. 
Fourth, \cref{sec:app-bivariate} provides an example where the Bayesian test's size may be strictly above $\alpha$ or equal to $\alpha$ depending on the sampling/posterior distribution. 
Another (less practically interesting) example is if $H_0 \colon \theta_2 \ne -\theta_1\textrm{ or }\theta_1=\theta_2=0$. 
Using the bivariate normal distribution above, if $\rho=-1$ and $\theta_1=\theta_2=0$, then the problem reduces to a point null hypothesis for which this Bayesian test has $100\%$ size. 
Conversely, if $\rho=1$, then the test's size is zero since the posterior probability of the complement of $\Theta_0$ is zero given any $(X_1,X_2)$.

For \cref{thm:super}:
the argument is the same as for \cref{thm:sub}, but with $\Theta_0$ larger than a half-space instead of smaller.
The first divergence from the proof of \cref{thm:sub} is the change from $\le$ and $<$ to $\ge$ and $>$ in
\begin{equation*}
\begin{split}
\Pr(\vecf{\theta} \in \Theta_0 \mid \phi_2(\vecf{X})=\vecf{p}^*) 
& \overbrace{\ge \Pr( \phi_2(\vecf{\theta}) \in \Phi_2 \mid \phi_2(\vecf{X})=\vecf{p}^*)}^{\textrm{since }\vecf{\theta} \in \Theta_0  \impliedby  \phi_2(\vecf{\theta}) \in \Phi_2}
\\&= \overbrace{\Pr(\phi_3(\phi_2(\vecf{\theta})) \le c_0 \mid \phi_2(\vecf{X})=\vecf{p}^*)}^{=\alpha\textrm{ by \cref{thm:half}}}
+ \overbrace{\Pr(\phi_2(\vecf{\theta}) \in \Delta \mid \phi_2(\vecf{X})=\vecf{p}^*)}^{>0\textrm{ by \cref{eqn:Delta-prob}}} 
\\& > \alpha . 
\end{split}
\end{equation*}
Analogous to before, there is an $\epsilon$-ball $\mathcal{B}$ around $\vecf{p}^*$ for which $\Pr(\vecf{\theta}\in\Theta_0 \mid \phi_2(\vecf{X})\in\mathcal{B})>\alpha$.
The same arguments as before lead to a statement parallel to \cref{eqn:ball-prob} (simply replacing $\mathcal{R}_a^\complement$ with $\mathcal{R}_a$),
\begin{equation}\label{eqn:ball-prob2}
\Pr\bigl( \phi_2(\vecf{X}) \in (\mathcal{B}\cap\mathcal{R}_a) \mid \vecf{\theta} \bigr) > 0\textrm{ for any }\vecf{\theta} . 
\end{equation}
Continuing in parallel, now $\vecf{\theta}\in\Theta_0 \impliedby \phi_2(\vecf{\theta})\in\Phi_2$ (instead of $\implies$), so
\begin{equation*}
\Pr(\phi_2(\vecf{\theta})\in\Phi_2 \mid \phi_2(\vecf{X})) \le \Pr(\vecf{\theta}\in\Theta_0 \mid \vecf{X}) ,
\end{equation*}
so \cref{eqn:Rb-implies-R2} is replaced by
\begin{equation}\label{eqn:Rb-impliedby-R2}
\phi_2(\vecf{X})\in\mathcal{R}_b
\impliedby  \vecf{X}\in\mathcal{R}_2 .
\end{equation}
The final argument becomes
\begin{align*}
\Pr(\vecf{X} \in \mathcal{R}_2 \mid \vecf{\theta}^*)
&\overbrace{%
\le \Pr(\phi_2(\vecf{X}) \in \mathcal{R}_b \mid \vecf{\theta}^*)
}^{\textrm{by \cref{eqn:Rb-impliedby-R2}}}
\\&\le \overbrace{\Pr(\phi_2(\vecf{X}) \in \mathcal{R}_a \mid \vecf{\theta}^*)}^{=\alpha\textrm{ by \cref{thm:half}}}
    - \overbrace{\Pr\bigl( \phi_2(\vecf{X}) \in (\mathcal{B}\cap\mathcal{R}_a) \mid \vecf{\theta}^* \bigr) }^{>0\textrm{ by \cref{eqn:ball-prob2}}}
\\&< \alpha . 
\qedhere
\end{align*}
\end{proof}

\section{{Bernstein--von Mises} theorems}
\label{sec:app-BvM}

A Bernstein--von Mises theorem establishes the asymptotic equivalence of the sampling and posterior distributions.
This can help verify \Cref{a:multi}, although it is not strictly necessary.
Results for posterior asymptotic normality date back to \citet{Laplace1820}, as cited in \citet[\S6.10, p.\ 515]{LehmannCasella1998}.

Parametric versions of the Bernstein--von Mises theorem are the oldest and can be found in textbooks.
They differ in regularity conditions and in how to quantify the distance between two distributions, but they share the requirement that the prior density be both continuous and positive at the true value.
For example, see Theorem 10.1 in \citet[\S10.2]{vanderVaart1998} and Theorems 20.1--3 in \citet[\S20.2]{DasGupta2008}, where Theorem 20.3 allows non-iid sampling.

General semiparametric versions of the Bernstein--von Mises theorem have been established, too.
For example, see \citet{Shen2002}, \citet{BickelKleijn2012}, and \citet{CastilloRousseau2015}, who allow non-iid sampling.
There are also papers for specific models like GMM and quantile regression; see 
\citet[Thm.\ G and footnote 13]{Hahn1997}, 
\citet[Thm.\ 2]{Kwan1999}, 
\citet[Prop.\ 1]{Kim2002}, 
\citet[Ex.\ 2]{Lancaster2003}, 
\citet[p.\ 36]{Schennach2005}, 
\citet[Sec.\ III.2]{Sims2010}, and 
\citet[Thm.\ 1]{Norets2015}, among others.

In the case of infinite-dimensional $\vecf{X}$ and $\vecf{\theta}$, as noted, an infinite-dimensional Bernstein--von Mises theorem is not necessary to satisfy \Cref{a:multi}, which only concerns the distribution of a scalar-valued functional.
However, the results and references below may be helpful or insightful in some cases.

For estimators of functions, it is common to have a (frequentist) Gaussian process limit with sample paths continuous with respect to the covariance semimetric; e.g., see \citet{vanderVaartWellner1996}.
A natural question is whether the (asymptotic, limit experiment) sampling and posterior distributions are ever equivalent in the sense of
\begin{equation*}
X(\cdot)-\theta(\cdot) \mid \theta(\cdot) \sim \mathbb{G} , \quad
\theta(\cdot)-X(\cdot) \mid X(\cdot) \sim \mathbb{G} , 
\end{equation*}
where $\mathbb{G}$ is a mean-zero Gaussian process with known covariance function.

Unfortunately, as discussed by \citet{Freedman1999} and others, such a Bernstein--von Mises result does not hold in great generality with infinite-dimensional spaces.
As explained by \citet[p.\ 1696]{HiranoPorter2009}, in finite dimensions the prior often behaves locally like Lebesgue measure (if its PDF is continuous and positive at the true parameter value), but in infinite-dimensional Banach spaces there is no analog of Lebesgue measure, let alone one that most priors would satisfy.

Fortunately, some nonparametric Bernstein--von Mises theorems do exist.
As in the (semi)parametric case, there are different ways to define ``asymptotically equivalent distributions''; see Definition 2 in \citet[p.\ 1950]{CastilloNickl2014} for an example.
The most general results to date seem to be provided by \citet{CastilloNickl2013,CastilloNickl2014}; see also Sections 12.4.1 and 12.2 of \citet{GhosalvanderVaart2017}.

One important special case where a Bernstein--von Mises theorem holds is for inference on a CDF.
On the frequentist side, assuming iid sampling,
\begin{equation}\label{eqn:inf-emp-process-limit}
\sqrt{n}\bigl( \hat{F}(\cdot) - F(\cdot) \bigr) 
\weaklyto   B\bigl(F(\cdot)\bigr) ,
\end{equation}
where $B(\cdot)$ is a standard Brownian bridge and $\weaklyto$ denotes weak convergence in $\ell^\infty(\bar\R)$; e.g., see \citet[Ex.\ 2.1.3]{vanderVaartWellner1996}.
For weak convergence under sequences $F_n(\cdot)\to F(\cdot)$, see Sections 2.8.3 and 3.11 and especially Theorem 3.10.12 in \citet{vanderVaartWellner1996}.
For a nonparametric Bayesian method using the Dirichlet process prior of \citet{Ferguson1973}, \citet[Thm.\ 2.1]{Lo1983} shows that a centered (at $\hat{F}(\cdot)$) and $\sqrt{n}$-scaled version of the posterior converges to the same limit as in \cref{eqn:inf-emp-process-limit} if the prior dominates $F(\cdot)$.
Even with an improper prior, i.e., using the Bayesian bootstrap of \citet{Rubin1981}, \citet[Thm.\ 2.1]{Lo1987} establishes the same result.
A closely related result is Theorem 12.2 of \citet{GhosalvanderVaart2017}.
An analogous conclusion is found in Theorem 4 of \citet{CastilloNickl2014}, but as a special case of their more general results.
They can provide Bernstein--von Mises theorems for (certain) collections of integral functionals of the PDF, $\int_0^1 g_t(x) f(x)\,\diff x$, where $f(\cdot)$ is the PDF with support $[0,1]$ and $t$ indexes the collection; for the CDF, $g_t(x)=\Ind{x \le t}$ for $t \in [0,1]$.

\section{Details on minimax risk}
\label{sec:app-minimax}

The following shows that an unbiased frequentist test with size $\alpha$ minimizes ``max'' (supremum) risk.
Let $\vecf{\theta} \in \Theta$, $H_0 \colon \vecf{\theta}\in\Theta_0\subset\Theta$, $H_1 \colon \vecf{\theta}\not\in\Theta_0$.  
From the frequentist perspective, given \cref{eqn:L}, the loss function given $\vecf{\theta}$ is $(1-\alpha)\Ind{\textrm{reject}}$ if $\vecf{\theta}\in\Theta_0$ or $\alpha\Ind{\textrm{accept}}$ if $\vecf{\theta}\not\in\Theta_0$.
Thus, risk (expected loss) given $\vecf{\theta}$ is $(1-\alpha)\Pr_{\vecf{\theta}}(\textrm{reject})$ if $\vecf{\theta}\in\Theta_0$ or $\alpha\Pr_{\vecf{\theta}}(\textrm{accept})$ if $\vecf{\theta}\not\in\Theta_0$, so ``minimax risk'' minimizes
\begin{equation}\label{eqn:max-risk}
\begin{split}
&\max\bigl\{ 
 (1-\alpha)\sup_{\vecf{\theta}\in\Theta_0}\Pr_{\vecf{\theta}}(\textrm{reject}), \;
 \alpha \sup_{\vecf{\theta}\not\in\Theta_0}\Pr_{\vecf{\theta}}(\textrm{accept})
    \bigr\} \\
&=
\max\bigl\{ 
 (1-\alpha)\sup_{\vecf{\theta}\in\Theta_0}\Pr_{\vecf{\theta}}(\textrm{reject}), \;
 \alpha [ 1 - \inf_{\vecf{\theta}\not\in\Theta_0}\Pr_{\vecf{\theta}}(\textrm{reject}) ]
    \bigr\}
,
\end{split}
\end{equation}
where $\Pr_{\vecf{\theta}}(\cdot)$ is probability under $\vecf{\theta}$. 
``Unbiased'' means the test satisfies
\begin{equation}\label{eqn:unbiased}
\sup_{\vecf{\theta}\in\Theta_0}\Pr_{\vecf{\theta}}(\textrm{reject})
 \le 
\inf_{\vecf{\theta}\not\in\Theta_0}\Pr_{\vecf{\theta}}(\textrm{reject}) .
\end{equation}
If the power function is continuous in $\vecf{\theta}$, then (writing $\partial\Theta_0$ for the boundary of $\Theta_0$, which is also the boundary of its complement) 
\begin{equation}\label{eqn:cts-power}
\sup_{\vecf{\theta}\in\Theta_0}
\Pr_{\vecf{\theta}}(\textrm{reject}) 
  \ge  
\sup_{\vecf{\theta}\in\partial\Theta_0}
\Pr_{\vecf{\theta}}(\textrm{reject}) 
  \ge  
\inf_{\vecf{\theta}\in\partial\Theta_0}
\Pr_{\vecf{\theta}}(\textrm{reject}) 
  \ge  
\inf_{\vecf{\theta}\not\in\Theta_0}
\Pr_{\vecf{\theta}}(\textrm{reject})
.
\end{equation}
Thus, combining \cref{eqn:unbiased,eqn:cts-power},
\begin{equation*}
\sup_{\vecf{\theta}\in\Theta_0}
\Pr_{\vecf{\theta}}(\textrm{reject}) 
=
\inf_{\vecf{\theta}\not\in\Theta_0}
\Pr_{\vecf{\theta}}(\textrm{reject}) 
.
\end{equation*}
Defining the test's size as
\begin{equation}\label{eqn:def-gamma0}
\gamma_0 \equiv \sup_{\vecf{\theta}\in\Theta_0}
\Pr_{\vecf{\theta}}(\textrm{reject}) ,
\end{equation}
\cref{eqn:max-risk} becomes
\begin{equation}
\max\{ (1-\alpha)\gamma_0, \alpha(1-\gamma_0) \} .
\end{equation}
The minimum possible value $\alpha(1-\alpha)$ is attained with size $\gamma_0=\alpha$: 
if instead $\gamma_0<\alpha$, then $\alpha(1-\gamma_0) > \alpha(1-\alpha)$, 
or if $\gamma_0>\alpha$, then $(1-\alpha)\gamma_0 > \alpha(1-\alpha)$.
Without unbiasedness, $\inf_{\vecf{\theta}\not\in\Theta_0}\Pr_{\vecf{\theta}}(\textrm{reject}) < \gamma_0$, so $\alpha [ 1 - \inf_{\vecf{\theta}\not\in\Theta_0}\Pr_{\vecf{\theta}}(\textrm{reject}) ]  >  \alpha(1-\gamma_0)$, increasing the minimax risk.
Thus, given \cref{eqn:cts-power}, a minimax risk decision rule is any unbiased frequentist hypothesis test with size $\alpha$.
See also \citet[Problem 1.10]{LehmannRomano2005text} on unbiased tests as minimax risk decision rules.

Without unbiasedness, the minimax-risk-optimal size of a test is above $\alpha$, but the magnitude of the difference is very small for conventional $\alpha$.
As a function of $\gamma_0$, let $\gamma_1(\gamma_0) \equiv \inf_{\vecf{\theta}\not\in\Theta_0}\Pr_{\vecf{\theta}}(\textrm{reject})$.
Then \cref{eqn:max-risk} equals
\begin{equation}\label{eqn:max-R0-R1}
\max\{ R_0 , R_1 \}
,\quad R_0 \equiv (1-\alpha)\gamma_0
,\quad R_1 \equiv \alpha(1-\gamma_1(\gamma_0)) .
\end{equation}
Continuity restricts $\gamma_1(\gamma_0) \le \gamma_0$, but we now drop the unbiasedness restriction $\gamma_1(\gamma_0) \ge \gamma_0$.
In the extreme, $\gamma_1(\gamma_0)=0$ for all $\gamma_0$, and the maximum risk is $\alpha$ for any test with $\gamma_0 \le \alpha/(1-\alpha)$, while maximum risk is larger than $\alpha$ if $\gamma_0>\alpha/(1-\alpha)$.
More commonly, $\gamma_1(\gamma_0)$ is strictly increasing in $\gamma_0$ but $\gamma_1(\alpha)<\alpha$.
This implies $R_0$ is increasing in $\gamma_0$ while $R_1$ is decreasing in $\gamma_0$, so the $\gamma_0$ that minimizes \cref{eqn:max-R0-R1} sets $R_0=R_1$.
If $\gamma_0=\alpha$, then $R_0<R_1$; if $\gamma_0=\alpha/(1-\alpha)$, then $R_0=\alpha>R_1$.
Hence, $R_0=R_1$ (and thus minimax risk) is achieved in between, at some $\gamma_0\in(\alpha,\alpha/(1-\alpha))$.
For example, rounding to two significant digits, if $\alpha=0.05$, then $\gamma_0\in(0.050,0.053)$, or if $\alpha=0.1$, then $\gamma_0\in(0.10,0.11)$. 
Such small divergence of $\gamma_0$ from $\alpha$ is almost imperceptible in practice.
Without unbiasedness, such a test is no longer a minimax rule (since the trivial randomized test ``reject with probability $\alpha$'' achieves lower $\alpha(1-\alpha)$ maximum risk), but the difference in the resulting maximum risk in \cref{eqn:max-R0-R1} is quantitatively small.

\section{Least favorable prior example}
\label{sec:app-LFP}

For intuition, consider a simple setting.
Let $X,\theta\in\R$, with $H_0 \colon \theta \le 0$ (so $\Theta_0=\{t:t\le0\}$), and $X\sim\NormDist(\theta,1)$, so the likelihood is $\ell(X \mid \theta)=(2\pi)^{-1/2}\exp\{-(X-\theta)^2/2\}$.
The prior sequence is $\pi_m(\theta)$ for $m=1,2,\ldots$, with posterior denoted $\pi_m(\cdot\mid X)$.
The corresponding Bayes rules are $\delta_m$, and the corresponding Bayes risks are $r_m$.
Given the limit $r_m\to r$, any decision rule whose supremum (over $\theta\in\R$) risk equals $r$ is minimax.
The posteriors are $\pi_m(\theta \mid X)$.
Below, the flat prior is shown not least favorable, and a least favorable sequence is derived.
The Bayes rules are from \cref{meth:Bayes}, i.e., the loss function still takes value $1-\alpha$ for type I error and $\alpha$ for type II error as in \cref{eqn:L}.

There is no single, fixed least favorable prior.
Consider the Bayes risk of the usual frequentist test that rejects when $X>\Phi^{-1}(1-\alpha)$, the $(1-\alpha)$-quantile of the $\NormDist(0,1)$ distribution.
Risk equals $\alpha(1-\alpha)$ if $\theta=0$, but risk is strictly below $\alpha(1-\alpha)$ for any $\theta\ne0$.
Thus, unless $\theta=0$ has prior probability one, the frequentist test's Bayes risk is strictly below $\alpha(1-\alpha)$, and by definition the Bayes rule's Bayes risk is weakly lower than that.
In the case where $\theta=0$ indeed has prior probability one, the posterior also puts all probability on $\theta=0$, so the Bayes rule is ``never reject,'' and the Bayes risk is zero (since ``never reject'' is always correct if $\theta=0$), also strictly below $\alpha(1-\alpha)$.

Next, consider the sequence of normal priors with improper limit, $\theta\sim\NormDist(0,m^2)$ as $m\to\infty$.
With large $m$, the posterior is approximately $\NormDist(X,1)$, so the Bayes rule rejects when $X>\Phi^{-1}(1-\alpha)$.
Again, when $\theta=0$, risk is $\alpha(1-\alpha)$, but when $\theta\ne0$, the risk is strictly below $\alpha(1-\alpha)$.
For example, when $\theta=-4$, $\Pr(X>\Phi^{-1}(1-\alpha))<0.00003$ for any $\alpha<0.5$.
Thus, at $\theta=-4$, risk is below $0.00003(1-\alpha)$.
The Bayes risk against $\pi_m$ averages the risk across all different values of $\theta$, so it will be below the maximum risk of $\alpha(1-\alpha)$.
Further, as $m\to\infty$, values of $\theta$ very far from zero become increasingly likely, so the Bayes risk becomes smaller and smaller, $r_m\to0$.
For example, even with $m=10$, $\Pr(\theta\le-4)=0.34$.
Coincidentally, the limiting (generalized) Bayes rule happens to be minimax, but this is not the least favorable prior sequence.
Arguably, it is the ``most favorable'' prior sequence because it places less and less weight on $\theta$ near the boundary point of $\Theta_0$ (i.e., near $\theta=0$), which is where risk is highest.

The least favorable prior sequence can be reverse-engineered from the unbiased frequentist test.
The prior must have $0<\pi_m(\Theta_0)<1$: both the null and alternative must have strictly positive probability, otherwise the test will always reject or always accept and have zero Bayes risk.
Given that the Bayes rule against $\pi_m$ rejects when $X>c_m$ for some $c_m\in\R$, there cannot be positive probability on both $\theta=0$ and $\theta<0$ (in the limit), or else the risk given $\theta<0$ is below risk given $\theta=0$.
The simplest remaining option is a two-point prior:
\begin{equation}\label{eqn:pi-m}
\pi_m(0) = p_m , \quad
\pi_m(a_m)=1-p_m\textrm{ for some }a_m>0 .
\end{equation}
In the limit, presumably the risk at $\theta=0$ is $\alpha(1-\alpha)$ because the rejection probability is $\alpha$.
If so, it must be that $a_m\downarrow0$; otherwise, the rejection probability at $a_m$ would remain strictly larger than $\alpha$, in which case risk at $\theta=a_m$ would remain strictly smaller than $\alpha(1-\alpha)$.

A sequence $p_m$ in \cref{eqn:pi-m} can be solved for given the sequence $a_m$.
For example, let $a_m=1/m\to0^+$.
We can solve for $p_m$ (below) such that the Bayes rule matches the frequentist test; it turns out $p_m\to\alpha$ as $m\to\infty$.
In that case, the risk at $\theta=0$ is $\alpha(1-\alpha)$ for all $m$, while the risk at $\theta=a_m$ approaches $\alpha(1-\alpha)$ from below as $m\to\infty$.
Since the Bayes risk averages these two risks (weighted by $p_m$ and $1-p_m$), and both approach $\alpha(1-\alpha)$, the limit of the sequence of Bayes risks is $\alpha(1-\alpha)$.
Thus, given that the usual frequentist test's supremum (over $\theta$) risk is $\alpha(1-\alpha)$ (the risk at $\theta=0$), the conclusion is that the frequentist test is a minimax risk decision rule.

The value of $p_m$ mentioned above can be found as follows.
The goal is for the Bayesian test to reject when $X>\Phi^{-1}(1-\alpha)$, same as the frequentist test.
Since the posterior probability of $\Theta_0$ decreases monotonically (and strictly) with $X$, it suffices to solve for $X$ such that the posterior equals $\alpha$ exactly, to find the boundary of the rejection region.
Since the prior has only two points with positive probability, so does the posterior.
Using the usual Bayes law of posterior proportional to likelihood times prior,
\begin{equation}\label{eqn:post-2pt-prop}
\pi_m(\theta=0 \mid X)   \propto b_{0m} \equiv p_m \phi(X) , \quad
\pi_m(\theta=a_m \mid X) \propto b_{1m} \equiv (1-p_m)\phi(X-a_m),
\end{equation}
where $\phi(\cdot)$ is the standard normal PDF.
Normalizing so these sum to one yields
\begin{equation}\label{eqn:post-2pt}
\pi_m(\theta=0 \mid X)   = b_{0m}/(b_{0m}+b_{1m}) , \quad
\pi_m(\theta=a_m \mid X) = b_{1m}/(b_{0m}+b_{1m}) .
\end{equation}
We want $\pi_m(\theta=0 \mid X)=\alpha$ when $X=\Phi^{-1}(1-\alpha)$, i.e.,
\begin{equation}
\alpha 
= \pi_m(0 \mid X=\Phi^{-1}(1-\alpha))
= \frac{p_m \phi\bigl(\Phi^{-1}(1-\alpha)\bigr)}%
       {p_m \phi\bigl(\Phi^{-1}(1-\alpha)\bigr)
        +(1-p_m)\phi\bigl(\Phi^{-1}(1-\alpha)-a_m\bigr)}
.
\end{equation}
This is easily solved numerically for $p_m$ given $\alpha$ and $a_m$.
Further, as $m\to\infty$ and $a_m\to0$, the right hand side converges to $p_m/(p_m+(1-p_m))=p_m$, so $p_m\downarrow\alpha$.
The risk at $\theta=0$ is $\alpha(1-\alpha)$, where $1-\alpha$ is the loss for type I error and $\alpha$ is the rejection probability, i.e., probability $X>\Phi^{-1}(1-\alpha)$ given $X\sim\NormDist(0,1)$.
The risk at $\theta=a_m$ is $\alpha\Phi(\Phi^{-1}(1-\alpha)-a_m)$, where $\alpha$ is the type II error loss and the rest is the probability of accepting $H_0$ (i.e., of $X\le\Phi^{-1}(1-\alpha)$) when $\theta=a_m$, given $X\sim\NormDist(a_m,1)$.
Altogether, the Bayes risk is
\begin{equation}\label{eqn:Bayes-risk-2pt}
 p_m (1-\alpha) \alpha
+(1-p_m) \alpha \Phi(\Phi^{-1}(1-\alpha)-a_m)
\to  (\alpha) (1-\alpha) \alpha
    +(1-\alpha) \alpha (1-\alpha)
= \alpha(1-\alpha) .
\end{equation}
Consequently, the test that rejects when $X>\Phi^{-1}(1-\alpha)$ is minimax since its supremum risk matches the Bayes risk limit in \cref{eqn:Bayes-risk-2pt}:
\begin{align*}
& \max\{ \sup_{\theta\le0} (1-\alpha)\Pr(\textrm{reject}\mid\theta) , \quad
       \sup_{\theta>0} \alpha\Pr(\textrm{accept}\mid\theta) \}
\\&=\max\{ \sup_{\theta\le0} (1-\alpha)\Pr(X>\Phi^{-1}(1-\alpha)\mid\theta) , \quad
\sup_{\theta>0} \alpha\Pr(X\le\Phi^{-1}(1-\alpha)\mid\theta) \}
\\&=\max\{ (1-\alpha)\Pr(X>\Phi^{-1}(1-\alpha)\mid\theta=0), \quad
 \alpha \lim_{\theta\downarrow0}\Pr(X\le\Phi^{-1}(1-\alpha)\mid\theta) \}
\\&=\max\{ (1-\alpha)\alpha , \quad \alpha(1-\alpha) \}
\\&=\alpha(1-\alpha) .
\end{align*}

\section{Example: bivariate normal sign equality test}
\label{sec:app-bivariate}

Here, we provide additional mathematical and simulation details for one of the examples mentioned in \cref{sec:results-discussion}. 
Consider a bivariate normal model where the null hypothesis is that the two parameters have the same sign (letting zero count either way), $H_0 \colon \theta_1 \theta_2 \ge 0$. 
Specifically,
\begin{equation}\label{eqn:bivariate-model}
\vecf{X} = (X_1,X_2)' \sim \NormDist(\vecf{\mu}, \matf{\Sigma}), \quad
\vecf{\mu} \equiv (\theta_1,\theta_2)', \quad
\matf{\Sigma} \equiv \left( \begin{array}{cc}1 & \rho \\ \rho & 1\end{array} \right)
.
\end{equation}
By symmetry, the test's size is the supremum of the test's rejection probability over $\theta_2 \in (-\infty,0]$ with $\theta_1=0$.
The parameter values leading to the biggest rejection probability depend on $\rho$; e.g., when $\rho=1$, the supremum comes from $\theta_2 \to -\infty$, whereas when $\rho=-1$, it comes at $\theta_2=0$.

To compute the Bayesian test, the posterior probability of $H_0$ given any $(X_1,X_2)$ must be computed.
Let the joint, conditional, and marginal PDFs of the $\NormDist(\vecf{0},\matf{\Sigma})$ distribution be
\begin{equation}\label{eqn:bivariate-PDFs}\begin{split}
\NormDist(\vecf{0},\matf{\Sigma}) : 
f(t_1, t_2) 
&\equiv \det(2\pi\matf{\Sigma})^{-1/2} \exp\left\{ -\frac{(t_1,t_2)\Sigma^{-1}(t_1,t_2)'}{2} \right\} \\
&= \frac{1}{2\pi \sqrt{1-\rho^2}} \exp\left\{ -\frac{t_1^2 + t_2^2 - 2 \rho t_1 t_2}{2(1-\rho^2)} \right\} , 
\\
\NormDist(\rho t_2, 1-\rho^2) : 
f(t_1 \mid t_2) &\equiv \frac{1}{\sqrt{2\pi(1-\rho^2)}} \exp\{ -(t_1-\rho t_2)^2 / (2-2\rho^2) \} , 
\\
\NormDist(0, 1) : 
f(t_2) &\equiv \frac{1}{\sqrt{2\pi}} \exp\{ -t_2^2/2 \} 
        = \phi(t_2) , 
\end{split}
\end{equation}
where $\phi(\cdot)$ is the standard normal PDF. 
Given \cref{eqn:bivariate-PDFs},
\begin{equation}
\int_{-\infty}^{-X_1} f(t_1 \mid t_2) \,\diff t_1
= \Phi\left( \frac{-X_1 - \rho t_2}{\sqrt{1-\rho^2}} \right) , 
\end{equation}
where $\Phi(\cdot)$ is the standard normal CDF. 
Using the prior expressions and the symmetry $f(t_1,t_2)=f(-t_1,-t_2)=f(t_2,t_1)$, 
\begin{align*}
\Pr( H_0 \mid X_1, X_2 )
&=\int_{-\infty}^{-X_2} 
  \int_{-\infty}^{-X_1} 
    f(t_1,t_2)\,\diff t_1 \,\diff t_2
+ \int_{-X_2}^{\infty} 
  \int_{-X_1}^{\infty} 
    f(t_1,t_2)\,\diff t_1 \,\diff t_2
\\&= 
  \int_{-\infty}^{-X_2} 
  \int_{-\infty}^{-X_1} 
    f(t_1 \mid t_2) f(t_2)\,\diff t_1 \,\diff t_2
+ \int_{-X_2}^{\infty} 
  \int_{-X_1}^{\infty} 
    f(t_1 \mid t_2) f(t_2)\,\diff t_1 \,\diff t_2
\\&= 
  \int_{-\infty}^{-X_2} 
    \Phi\left( \frac{-X_1 - \rho t_2}{\sqrt{1-\rho^2}} \right)
    \phi(t_2)\,\diff t_2
+ \int_{-X_2}^{\infty} 
    \left[ 1 - \Phi\left( \frac{-X_1 - \rho t_2}{\sqrt{1-\rho^2}} \right) \right]
    \phi(t_2)\,\diff t_2
.
\end{align*}
Unfortunately, there is no closed-form expression for this $\Pr( H_0 \mid X_1, X_2 )$. 
However, it is easily simulated. 
The function $\Phi(\cdot)$ is available in any modern statistical software (e.g., R). 
After drawing many $Z_j  \iid  \NormDist(0,1)$ for $j=1,\ldots,J$, letting 
$W_j \equiv \Phi((-X_1-\rho Z)/\sqrt{1-\rho^2})$ if $Z<-X_2$ and $W_j \equiv 1 - \Phi((-X_1-\rho Z)/\sqrt{1-\rho^2})$ if $Z \ge -X_2$, then $\Pr( H_0 \mid X_1, X_2 ) \approx J^{-1} \sum_{j=1}^{J} W_j$, with the approximation error going to zero as $J\to\infty$.

\begin{table}[htbp]
\centering
\caption{\label{tab:sim-bivariate-RP}Bayesian test type I error rates, $H_0 \colon \theta_1 \theta_2 \ge 0$, $\theta_1=0$, $\alpha=0.1$, $\numnornd{10000}$ replications.}
\sisetup{round-precision=2}
\begin{tabular}[c]{S[round-precision=2,round-mode=places] S[table-format=1.3,round-precision=3,round-mode=places]S[table-format=1.3,round-precision=3,round-mode=places]S[table-format=1.3,round-precision=3,round-mode=places]S[table-format=1.3,round-precision=3,round-mode=places]S[table-format=1.3,round-precision=3,round-mode=places]S[table-format=1.3,round-precision=3,round-mode=places]S[table-format=1.3,round-precision=3,round-mode=places]}
\toprule
 & \multicolumn{7}{c}{$\rho$} \\
\cmidrule{2-8}
{$\theta_2$} & {$-1$} & {$-0.99$} & {$-0.9$} & {$-0.5$} & {$0$} & {$0.5$} & {$1$} \\
\midrule
  0.00 & 1.00000 & 0.96640 & 0.26390 & 0.05760 & 0.01050 & 0.00010 & 0.00000 \\
 -0.25 & 1.00000 & 0.71740 & 0.25260 & 0.05860 & 0.01130 & 0.00010 & 0.00000 \\
 -0.50 & 0.25315 & 0.30480 & 0.22120 & 0.06210 & 0.01430 & 0.00060 & 0.00000 \\
 -4.00 & 0.10000 & 0.10180 & 0.10150 & 0.10090 & 0.09920 & 0.09400 & 0.09303 \\
-10.00 & 0.10000 & 0.10180 & 0.10150 & 0.10080 & 0.10020 & 0.10020 & 0.10000 \\
\midrule
{Max}  & 1.00000 & 0.96640 & 0.26390 & 0.10090 & 0.10020 & 0.10020 & 0.10000 \\
\bottomrule
\end{tabular}
\end{table}

\Cref{tab:sim-bivariate-RP} shows type I error rates of the Bayesian test of $H_0 \colon \theta_1 \theta_2 \ge 0$ for different $\rho$ and different $\theta_2$, with $\theta_1=0$. 
When $\rho=-1$, it reduces to a one-dimensional setting with a finite, convex $\Theta_0$, as in \cref{thm:sub}; size equals one, well above $\alpha$. 
When $\rho=1$, it is also essentially one-dimensional, but with non-convex $\Theta_0$, as in \cref{thm:non}; the type I error rate is (near) zero at many points on the boundary of $\Theta_0$, but size equals $\alpha$. 
When $\rho\in(-1,1)$, we are also in the setting of \cref{thm:non}; size is strictly above $\alpha$ when $\rho$ is close to $-1$, but decreases to $\alpha$ somewhere between $\rho=-0.9$ and $\rho=-0.5$.

\section{Derivation of translog constraints}
\label{sec:app-translog}

The Hessian is a nonlinear function of the translog parameters, and it depends on $(y,\vecf{w})$. 
Letting\footnote{Some notation is from \citet{ODonnellCoelli2005}.}  
\begin{equation}\label{eqn:r}
r_k \equiv \pD{\ln(C)}{\ln(w_k)}
= a_{yk}\ln(y) + b_k + \sum_{j=1}^{3}b_{jk}\ln(w_j) , 
\end{equation}
a general element of $\matf{H}$ is 
\begin{align*}
H_{mk} 
&= \frac{\partial^2 C}{\partial w_m \partial w_k}
 = \pD{}{w_m} \pD{C}{w_k}
 = \pD{}{w_m} ( r_k C / w_k )
 = \pD{r_k}{w_m} (C/w_k) 
    +\pD{C}{w_m} (r_k/w_k)
    +\pD{w_k^{-1}}{w_m} (r_k C)
\\&= (b_{mk}/w_m) (C/w_k) 
    +r_m (C/w_m) (r_k/w_k)
    -\Ind{k=m}w_k^{-2} (r_k C)
\\&= C \frac{b_{mk} + r_m r_k - \Ind{k=m}r_k}{w_m w_k} 
. 
\end{align*}
Since each element is proportional to $C>0$, the value of $C$ does not affect whether or not $\matf{H}$ is NSD: $\matf{H}$ is NSD iff $\matf{H}/C$ is NSD. 
This may be helpful if the translog parameters are estimated from cost share equations and $C$ is not directly observed. 

The local NSD condition in \cref{eqn:H0-translog-local} corresponds to a set of parameter values much smaller than a half-space. 
A necessary (but not sufficient) condition for NSD is that all the principal minors of order $p=1$ are non-positive, i.e., that $H_{11}\le0$, $H_{22}\le0$, and $H_{33}\le0$. 
In terms of the parameters, using \cref{eqn:r}, $H_{11}\le0$ iff $b_{11}+r_1^2-r_1\le0$, i.e., $b_{11}\le r_1(1-r_1)$. 
With $(y,\vecf{w})=(1,1,1,1)$, $r_k=b_k$, so $H_{11}\le0$ iff $b_{11}\le b_1(1-b_1)$. 
After imposing symmetry ($b_{mk}=b_{km}$) and homogeneity of degree one in input prices ($b_{m1}+b_{m2}+b_{m3}=0$, $m=1,2,3$), all $b_{mk}$ can be written in terms of $b_{11}$, $b_{12}$, and $b_{22}$: $b_{21}=b_{12}$, $b_{13} = -b_{11} -b_{12}$, etc. 
Also from homogeneity, $b_1+b_2+b_3=1$, and from monotonicity, $b_k=r_k\ge0$, so $0\le b_1\le1$. 
Thus, $b_1(1-b_1)\in[0,0.25]$, so $b_{11}\le b_1(1-b_1)$ is larger than the half-space defined by $b_{11}\le0$ but smaller than the half-space defined by $b_{11}\le0.25$. 
A similar argument for $H_{22}\le0$ at $(1,1,1,1)$ yields $b_{22}\le b_2(1-b_2)\le0.25$. 
From the constraints on $H_{11}$ and $H_{22}$ alone, $\Theta_0$ is a subset of the ``quarter-space'' defined by $b_{11}\le0.25$ and $b_{22}\le0.25$. 
In the notation of \cref{thm:1}, we could use $\phi(\vecf{\theta})=b_{11}$ (or $b_{22}$) and $c_0=0.25$. 
Adding the constraints for the other principal minors of $\matf{H}$ makes $\Theta_0$ even smaller. 

Since the local concavity $H_0$ in \cref{eqn:H0-translog-local} corresponds to a subset of a \emph{quarter}-space in the parameter space, \cref{thm:sub} suggests that we expect the Bayesian test's size to exceed $\alpha$. 
The results in \cref{tab:sim-translog} show this to be the case here.

\section{Example: first-order stochastic dominance}
\label{sec:app-SD1}

For testing first-order stochastic dominance (SD1), let $X_i \iid F_X(\cdot)$, $Y_i \iid F_Y(\cdot)$ and independent of the $X_i$ sample, and $F_0(\cdot)$ is non-random, where all distributions are continuous. 
One-sample SD1 is $F_X(\cdot)\le F_0(\cdot)$; two-sample SD1 is $F_X(\cdot)\le F_Y(\cdot)$; ``non-SD1'' means SD1 is not satisfied.

First, we show how \cref{thm:1} applies to SD1. 
Second, we provide analytic results from the limit experiment. 
Third, we show simulated finite-sample results.

\subsection{SD1: application of \texorpdfstring{\cref{thm:1}}{Theorem \ref{thm:1}}}
\label{sec:ex-SD-thm}

As we show below, \cref{thm:1} implies that the Bayesian test's asymptotic size is strictly above $\alpha$ when the null hypothesis is SD1 but that this may not hold when the null is non-SD1. 
This subsection shows how SD1 and non-SD1 satisfy the conditions of \cref{thm:sub} and \cref{thm:non}, respectively.

Although $F_X(\cdot)$ is infinite-dimensional, only finite-dimensional marginal distributions are required to apply \cref{thm:1}. 
Consider the simpler \cref{thm:subeq} first. 
Let $X_n(\cdot) \equiv \sqrt{n}\bigl(\hat{F}_X(\cdot)-F_{0,n}(\cdot)\bigr)$ and $\theta_n(\cdot) \equiv \sqrt{n}\bigl(F_X(\cdot)-F_{0,n}(\cdot)\bigr)$, where (from the frequentist view) $\theta_n(\cdot) \to \theta(\cdot)$, the local mean parameter. 
SD1 is equivalent to $\theta(\cdot) \le 0(\cdot)$. 
Although limits of the full infinite-dimensional sequences are tractable (see \cref{sec:ex-SD-limit}), only a scalar-valued functional is needed for \cref{thm:subeq}. 
Let $\phi(\theta(\cdot))=\theta(r)$ for some (any) $r \in \R$, so $\theta(\cdot) \le 0(\cdot)  \implies  \phi(\theta(\cdot)) \le 0$, satisfying the condition of \cref{thm:subeq} that $\Theta_0 \subseteq \{ \theta(\cdot) : \phi(\theta(\cdot)) \le 0 \}$. 
Let $D_n \equiv \phi(X_n(\cdot)) = \sqrt{n}(\hat{F}_X(r)-F_{0,n}(r))$ and $\gamma_n \equiv \phi(\theta_n(\cdot)) = \sqrt{n}(F_X(r)-F_{0,n}(r))$. 
Writing $Z_i=\Ind{ X_i \le r }$, then $F_X(r)=\E(Z_i)$ and $\hat{F}_X(r)=n^{-1}\sum_{i=1}^{n} Z_i$, i.e., we are concerned only with the mean of a random variable. 
The asymptotic sampling distribution of $D_n-\gamma_n$ is $\NormDist\bigl(0,F_X(r)[1-F_X(r)]\bigr)$, satisfying the continuity, support, and symmetry conditions in \Cref{a:multi}. 
The remainder of \cref{a:multi} is satisfied if a semiparametric Bernstein--von Mises theorem for the mean holds. 
This and even stronger results hold with a Dirichlet process prior, as in \citet{Lo1983}.

For \cref{thm:sub}, we only need strengthen the Bernstein--von Mises theorem from a scalar to a bivariate vector ($d=2$), which again holds with a Dirichlet process prior \citep{Lo1983}, for example. 
Here, let $\phi_2(\theta(\cdot)) = (\theta(r),\theta(t))$ for $r,t\in\R$, and $\phi_3(\vecf{p})=p_1$ so $\phi_3(\phi_2(\theta(\cdot))) = \theta(r) = \phi(\theta(\cdot))$. 
Also, $\theta(\cdot) \le 0(\cdot)  \implies  (\theta(r), \theta(t)) \le \vecf{0}$, satisfying the condition in \cref{thm:sub} on $\Phi_2 = \{ \vecf{p} : p_1 \le 0, p_2 \le 0 \}$. 
Continuing to follow the notation from \cref{thm:sub}, $\Delta = \{\vecf{p} : p_1 \le 0, p_2 > 0 \}$, which has positive (indeed infinite) Lebesgue measure as required. 
The bivariate asymptotic distribution $F_2$ is bivariate normal, which again satisfies \cref{a:multi} as well as the continuity and strictly positive PDF requirement.

For testing non-SD1, \cref{thm:non} applies. 
Non-SD1 is satisfied in the entire half-space $\{ \theta(\cdot) : \theta(r) \ge 0 \}$, as well as most of the complement half-space (e.g., if $\theta(r)<0$ but $\theta(t) \ge 0$), so it cannot be contained in any half-space.

For the two-sample setting, the infinite-dimensional limits in \cref{sec:ex-SD-limit} are more than sufficient to establish the scalar and bivariate conditions of \cref{thm:sub}, and again non-SD1 cannot be contained in any half-space.

\subsection{SD1: results from limit experiment}
\label{sec:ex-SD-limit}

We derive the infinite-dimensional limit experiment and then compute certain results.
Continuing some notation from \cref{sec:ex-SD-thm}, consider the one-sample setup with $X_i \iid F_X(\cdot)$.
Again let $\theta_n(\cdot) \equiv \sqrt{n}\bigl(F_X(\cdot)-F_{0,n}(\cdot)\bigr) \to \theta(\cdot)$, the local parameter.
SD1 of $F_X$ over $F_0$ can be written equivalently as
\begin{equation}\label{eqn:SD1}
F_X \SD{1} F_0  \iff  F_X(\cdot) \le F_0(\cdot)  \iff \theta(\cdot) \le 0(\cdot) . 
\end{equation}
Since (by Donsker's theorem) 
$\sqrt{n}\bigl(\hat{F}_X(\cdot)-F_X(\cdot)\bigr) \weaklyto B\bigl(F_X(\cdot)\bigr)$ 
for standard Brownian bridge $B(\cdot)$, similar to \cref{eqn:gen-drift-mu0},
\begin{align}\notag
X_n(\cdot) &\equiv \sqrt{n}\bigl(\hat{F}_X(\cdot)-F_{0,n}(\cdot)\bigr)
\\&= \sqrt{n}\bigl(\hat{F}_X(\cdot)-F_X(\cdot)\bigr) 
 +\sqrt{n}\bigl(F_X(\cdot)-F_{0,n}(\cdot)\bigr)
 \notag
\\&\weaklyto B\bigl(F_X(\cdot)\bigr)+\theta(\cdot) , 
\label{eqn:SD1-X-limit}
\end{align}
so the limit experiment has 
$X(\cdot)-\theta(\cdot) \mid \theta(\cdot) \sim B\bigl(F_X(\cdot)\bigr)$.
Note $B\bigl(F_X(\cdot)\bigr)$ is a mean-zero Gaussian process with covariance function $\Cov(t_1,t_2)=F_X(t_1)[1-F_X(t_2)]$ for $t_1\le t_2$.
Although $F_X(\cdot)$ is unknown, $\hat{F}_X(\cdot) \stackrel{a.s.}{\to} F_X(\cdot)$ uniformly by the Glivenko--Cantelli theorem, so asymptotically the covariance function is known while the (local) mean function $\theta(\cdot)$ remains unknown.
Analogously, for the posterior, similar to \cref{eqn:gen-drift-mu0-posterior},
\begin{equation}
\theta_n(\cdot) - X_n(\cdot) 
= \sqrt{n}\bigl(F_X(\cdot)-\hat{F}_X(\cdot)\bigr) 
\weaklyto B\bigl(F_X(\cdot)\bigr) , 
\label{eqn:SD1-theta-limit}
\end{equation}
using the Bernstein--von Mises theorem in \citet{Lo1983,Lo1987} for Dirichlet process prior Bayesian inference or Theorem 4 of \citet{CastilloNickl2014}.

The two-sample setting is similar since we assume the samples are independent.
For notational simplicity, assume both samples have $n$ observations.
Let $\Delta(\cdot) \equiv F_X(\cdot) - F_Y(\cdot)$, the true CDF difference function.
Let $\Delta_{0,n}(\cdot)$ be the centering functions satisfying $\sqrt{n}(\Delta(\cdot) - \Delta_{0,n}(\cdot)) \to \theta(\cdot)$, the local parameter.
SD1 of $F_X$ over $F_Y$ is
\begin{equation}\label{eqn:SD1-2s}
F_X \SD{1} F_Y  \iff  F_X(\cdot) \le F_Y(\cdot)  \iff  \theta(\cdot) \le 0(\cdot) . 
\end{equation}
For the sampling distribution,
\begin{align}\notag
X_n(\cdot)
&\equiv \sqrt{n}( \hat{F}_X(\cdot) - \hat{F}_Y(\cdot) - \Delta_{0,n}(\cdot) )
\\&= \sqrt{n}(\hat{F}_X(\cdot) - F_X(\cdot))
 -\sqrt{n}(\hat{F}_Y(\cdot) - F_Y(\cdot))
 +\sqrt{n}(\Delta(\cdot) - \Delta_{0,n}(\cdot))
\notag
\\&
\weaklyto B_1(F_X(\cdot))
         -B_2(F_Y(\cdot))
         +\theta(\cdot) 
, 
\label{eqn:SD1-X-limit-2s}
\end{align}
where $B_1(\cdot)$ and $B_2(\cdot)$ are independent standard Brownian bridges.
For the posterior, using the independence of samples and Bernstein--von Mises theorem,
\begin{align}\notag
\theta_n(\cdot) - X_n(\cdot) 
&= \sqrt{n}(F_X(\cdot) - \hat{F}_X(\cdot))
  -\sqrt{n}(F_Y(\cdot) - \hat{F}_Y(\cdot))
\\&
\weaklyto B_1(F_X(\cdot))
         -B_2(F_Y(\cdot))
. \label{eqn:SD1-theta-limit-2s}
\end{align}

Three results are now stated, followed by proofs.

First, consider the Bayesian posterior probability of SD1 when $X(\cdot)=0(\cdot)$. 
The finite-sample analog is when $\hat{F}_X(\cdot) \approx F_0(\cdot)$ or $\hat{F}_X(\cdot) \approx \hat{F}_Y(\cdot)$. 
The value $0(\cdot)$ is at the very ``corner'' of $\Theta_0$, and it is a very pointy corner, so a ball centered at $0(\cdot)$ contains very little of $\Theta_0$. 
In fact, ``very little'' means ``zero probability,'' as the next result states.

\begin{proposition}\label{prop:SD1-PrH0}
Consider the limit experiment posterior for one-sample SD1 testing in \cref{eqn:SD1-theta-limit} and for two-sample SD1 testing in \cref{eqn:SD1-theta-limit-2s}. 
Given $X(\cdot)=0(\cdot)$, the posterior probability of SD1 is zero in both one-sample and two-sample settings. 
\end{proposition}

Second, similar intuition and arguments lead to the SD1 Bayesian test's size being one. 
\begin{proposition}\label{prop:SD1-size}
Consider limit experiments for one-sample SD1 testing in \cref{eqn:SD1-X-limit,eqn:SD1-theta-limit} and for two-sample SD1 testing in \cref{eqn:SD1-X-limit-2s,eqn:SD1-theta-limit-2s}. 
Consider the Bayesian test from \cref{meth:Bayes} that rejects $H_0 \colon \theta(\cdot) \le 0(\cdot)$ iff the posterior probability of $H_0$ is below $\alpha$. 
Then, the Bayesian test's size equals one, with type I error rate equal to one when $\theta(\cdot)=0(\cdot)$. 
\end{proposition}

Third, the following result for non-SD1 rejection probability is immediate. 
\begin{corollary}\label{cor:nonSD1-RP}
Consider same setup as in \cref{prop:SD1-size}. 
When $\theta(\cdot)=0(\cdot)$, the Bayesian test's probability of rejecting non-SD1 is zero. 
\end{corollary}

\begin{proof}[\bfseries Proof of \cref{prop:SD1-PrH0}]
In the one-sample case, using \cref{eqn:SD1-theta-limit},
\begin{align}\notag
     \Pr( F_X \SD{1} F_0 \mid X(\cdot)=0(\cdot) ) 
  &= \Pr( \theta(\cdot) \le 0(\cdot)  \mid  X(\cdot)=0(\cdot) ) 
   = \Pr( B(F_X(\cdot)) \le 0(\cdot) )
\\&= \Pr( B(\cdot) \le 0(\cdot) )
   = \Pr( \sup_{t \in [0,1]} B(t) \le 0 )
   = 0 . 
\label{eqn:Pr-sup-B-neg}
\end{align}
The final equality holds because the distribution of the supremum of a mean-zero Brownian bridge is continuous and has non-negative support; e.g., see Theorem 2 in \citet{Smirnov1939b} or equation (1.1) in \citet{BirnbaumTingey1951}.

For two-sample testing, the result in \cref{eqn:Pr-sup-B-neg} extends readily if we assume $F_X=F_Y$ since then $B_1(F_X(\cdot))-B_2(F_Y(\cdot))=B_1(F(\cdot))-B_2(F(\cdot)) \stackrel{d}{=} \sqrt{2}B(F(\cdot))$ for another independent Brownian bridge $B(\cdot)$.
More generally,%
\footnote{Thanks to Iosif Pinelis for help extending to $F_X \ne F_Y$: \url{https://mathoverflow.net/a/292716/120669}}
let
\begin{equation*}
T(\cdot) \equiv B_1(F_X(\cdot)) - B_2(F_Y(\cdot)) ,
\end{equation*}
the distribution of $\theta(\cdot)$ conditional on $X(\cdot)=0(\cdot)$. 
Using \cref{eqn:SD1-theta-limit-2s},
\begin{align}\notag
     \Pr( F_X \SD{1} F_Y \mid X(\cdot)=0(\cdot) ) 
  &= \Pr( \theta(\cdot) \le 0(\cdot)  \mid  X(\cdot)=0(\cdot) ) 
   = \Pr( T(\cdot) \le 0(\cdot) )
\\&= \Pr( \sup_{r\in\R} T(r) \le 0 )
.
\label{eqn:PrH0-2s-Xeq0}
\end{align}
Let $L$ denote the smaller of the lower bounds of the distributions $F_X$ and $F_Y$, allowing $L=-\infty$ if both have unbounded support.
Let $W(\cdot)$, $W_1(\cdot)$, and $W_2(\cdot)$ denote independent standard Brownian motion processes.
We may write
\begin{gather*}
B_1(t) = W_1(t) - t W_1(1), \quad 
B_2(t) = W_2(t) - t W_2(1), \\
V(\cdot) \equiv W_1(F_X(\cdot)) - W_2(F_Y(\cdot)) 
  \stackrel{d}{=} \sqrt{2} W((F_X(\cdot)+F_Y(\cdot))/2), \\
Z(\cdot) \equiv V(\cdot) - T(\cdot)
  = F_X(\cdot)W_1(1) - F_Y(\cdot)W_2(1) . 
\end{gather*}
Looking at $T(r)=V(r)-Z(r)$ as $r \downarrow L$, the $Z(r)$ becomes negligibly small, while the $V(r)$ varies sufficiently to attain a strictly positive supremum almost surely.
Specifically,
\begin{equation*}
\lim_{r \downarrow L} \frac{Z(r)}{\sqrt{F_X(r)+F_Y(r)}} = 0 , 
\end{equation*}
so continuing from \cref{eqn:PrH0-2s-Xeq0} with the $\le$ changed to $>$,
\begin{equation*}
\begin{split}
\Pr( \sup_{r\in\R} T(r) > 0 )
  &\ge \Pr\left( \limsup_{r \downarrow L} \frac{T(r)}{\sqrt{F_X(r)+F_Y(r)}} = \infty \right)
\\&= \Pr\left( \limsup_{r \downarrow L} \frac{V(r)-Z(r)}{\sqrt{F_X(r)+F_Y(r)}} = \infty \right)
\\&= \Pr\left( \limsup_{r \downarrow L} \frac{V(r)}{\sqrt{F_X(r)+F_Y(r)}} = \infty \right)
\\&= \Pr\left( \limsup_{r \downarrow L} \frac{\sqrt{2}}{\sqrt{F_X(r)+F_Y(r)}} W\left(\frac{F_X(\cdot)+F_Y(\cdot)}{2}\right) = \infty \right)
= 1
\end{split}
\end{equation*}
by the (local) law of iterated logarithm.\footnote{E.g., Corollary 5.3 in \url{https://www.stat.berkeley.edu/~peres/bmbook.pdf}}
\end{proof}

\begin{proof}[\bfseries Proof of \cref{prop:SD1-size}]
We show that the type I error rate is one when $\theta(\cdot)=0$, which directly implies the size is one, too. 
For the two-sample case, $\theta(\cdot)=0(\cdot)$ implies $F_X(\cdot)=F_Y(\cdot)$, so the limit experiment simplifies since $B_1(F_X(\cdot))-B_2(F_Y(\cdot)) \stackrel{d}{=} \sqrt{2} B(F(\cdot))$ for $F(\cdot)=F_X(\cdot)=F_Y(\cdot)$ and standard Brownian bridge $B(\cdot)$. 
Thus, both one-sample and two-sample limiting distributions can be written as $cB(\cdot)$, where $c=1$ for one-sample and $c=\sqrt{2}$ for two-sample.

The Bayesian test rejects when the posterior is below $\alpha$, so the probability of \emph{not} rejecting when $\theta(\cdot)=0$ is
\begin{align*}
&
\Pr( X(\cdot) \in \{ x(\cdot) : \Pr( \theta(\cdot) \le 0(\cdot) \mid X(\cdot)=x(\cdot) ) > \alpha \} \mid \theta(\cdot)=0(\cdot) ) 
\\&= 
\Pr( X(\cdot) \in \{ x(\cdot) : \Pr( c B(F(\cdot)) \le -x(\cdot) ) > \alpha \} \mid \theta(\cdot)=0(\cdot) ) . 
\end{align*}
This can be shown to be zero via the unconditional probability
\begin{equation}\label{eqn:Pr-BB-sum-neg}
\Pr( c B_1(F(\cdot)) + c B_2(F(\cdot)) \le 0(\cdot) )
= \Pr( \sqrt{2}c B(F(\cdot)) \le 0(\cdot) )
= \Pr( B(\cdot) \le 0(\cdot) )
= 0 ,
\end{equation}
again using \cref{eqn:Pr-sup-B-neg}, where $B_1(\cdot)$, $B_2(\cdot)$, and $B(\cdot)$ are independent standard Brownian bridges.
For any set $S$, with $X(\cdot) \sim c B_2(F(\cdot))$ since $\theta(\cdot)=0(\cdot)$,
\begin{equation*}\begin{split}
&
\Pr( \theta(\cdot) \le 0(\cdot) \mid X(\cdot) \in S)
\Pr( X(\cdot) \in S) 
\\&=
\Pr( c B_1(F(\cdot)) + X(\cdot) \le 0(\cdot) \mid X(\cdot) \in S)
\Pr( X(\cdot) \in S ) 
\\&=
\Pr( c B_1(F(\cdot)) + X(\cdot) \le 0(\cdot)\textrm{ and }X(\cdot) \in S )
\\&\le
\Pr( c B_1(F(\cdot)) + X(\cdot) \le 0(\cdot) )
\\&=
\Pr( c B_1(F(\cdot)) + c B_2(F(\cdot)) \le 0(\cdot) )
\\&= 0
\end{split}
\end{equation*}
by \cref{eqn:Pr-BB-sum-neg}, where $c B_2(F(\cdot))$ is the sampling distribution of $X(\cdot)$ given $\theta(\cdot)=0(\cdot)$, and in the posterior $\theta(\cdot) \sim c B_1(F(\cdot)) + X(\cdot)$.
Consequently,
\begin{equation}\label{eqn:Pr-theta-neg-uncond}
\Pr( \theta(\cdot) \le 0(\cdot) \mid X(\cdot) \in S)
\Pr( X(\cdot) \in S) 
\le 0 . 
\end{equation}
Specifically, let $S$ be the complement of the test's rejection region:
\begin{equation*}
S = \{ x(\cdot) : \Pr( \theta(\cdot) \le 0(\cdot) \mid X(\cdot)=x(\cdot) ) > \alpha \} . 
\end{equation*}
If $\Pr( X(\cdot) \in S)>0$, then the left-hand side of \cref{eqn:Pr-theta-neg-uncond} is the product of two strictly positive terms (assuming $\alpha>0$), which is strictly positive.
This contradicts \cref{eqn:Pr-theta-neg-uncond} since the right-hand side is zero. 
Consequently, $\Pr( X(\cdot) \in S)=0$ and thus $\Pr( X(\cdot) \not\in S)=1$, i.e., the rejection probability is one.
(This does not mean $S$ is empty, just that it is a zero-probability set.)
\end{proof}

\begin{proof}[\bfseries Proof of \cref{cor:nonSD1-RP}]
With $\theta(\cdot)=0(\cdot)$, by \cref{prop:SD1-size}, the probability of rejecting SD1 is $100\%$; that is, the posterior probability of SD1 is below $\alpha$ with probability one (with respect to the distribution of $X(\cdot)$).
Since the posterior probabilities of SD1 and non-SD1 sum to one, this implies the posterior probability of non-SD1 is above $1-\alpha$ (and thus the test does not reject) with probability one.
\end{proof}

\subsection{SD1: finite-sample simulations}
\label{sec:ex-SD-sims}

The following simulation results reflect the theoretical results from the limit experiment (discussed in \cref{sec:ex-SD-thm}). 
Code for replication is provided.

\Cref{tab:sim-pval} shows Bayesian posterior probabilities of SD1 and non-SD1 in datasets near the ``corner'' of SD1, similar to the setup of \cref{prop:SD1-PrH0}.
In the one-sample case, this means $\hat{F}_X(\cdot)$ nearly equals the $\UnifDist(0,1)$ CDF. 
In the two-sample case, this means $\hat{F}_X(\cdot)$ nearly equals $\hat{F}_Y(\cdot)$. 
Specifically, we set $X_i=i/(n+1)$ for $i=1,\ldots,n$, and in the two-sample case, $Y_i=i/n$ for $i=1,\ldots,n-1$ (there are $n-1$ observations in the second sample). 
When $n<\infty$, the Bayesian bootstrap variant of \citet{Banks1988} is used. 
When $n=\infty$, the results are from \cref{prop:SD1-PrH0}. 

\begin{table}[htbp]
\centering
\caption{\label{tab:sim-pval}Bayesian posterior probabilities of $H_0 \colon X \SD{1} \UnifDist(0,1)$ and $H_0 \colon X \SD{1} Y$.} 
\begin{tabular}{rrrr}
\toprule
&     &  \multicolumn{2}{c}{Comparison distribution} \\
\cmidrule{3-4}
\multicolumn{1}{r}{$H_0$} & 
\multicolumn{1}{r}{$n$}   &
$\UnifDist(0,1)$ & 
$Y$ \\
\midrule
    SD1 &     $10$ & 0.103 & 0.097 \\
    SD1 &     $40$ & 0.028 & 0.025 \\
    SD1 &    $100$ & 0.009 & 0.010 \\
    SD1 & $\infty$ & 0.000 & 0.000 \\[2pt]
non-SD1 &     $10$ & 0.897 & 0.903 \\
non-SD1 &     $40$ & 0.972 & 0.975 \\
non-SD1 &    $100$ & 0.991 & 0.990 \\
non-SD1 & $\infty$ & 1.000 & 1.000 \\
\bottomrule
\end{tabular}
\end{table}

\Cref{tab:sim-pval} illustrates the Bayesian interpretation of a draw near $X(\cdot)=0(\cdot)$. 
This interpretation differs greatly from a frequentist interpretation and illuminates the rejection probabilities seen in \cref{tab:sim-RP}. 
The same intuition before \cref{prop:SD1-PrH0} applies here. 
Consequently, when $\hat{F}_X(\cdot) \approx F_0(\cdot)$, or when $\hat{F}_X(\cdot)\approx\hat{F}_Y(\cdot)$, the Bayesian posterior places nearly zero probability on SD1 and (equivalently) almost all probability on non-SD1.  
\Cref{tab:sim-pval} shows finite-sample posterior probabilities when $n=100$ to be very close to the limit as $n \to \infty$. 
Opposite the Bayesian interpretation, a frequentist $p$-value for the null of SD1 would be near one when the estimated $\hat{F}_X(\cdot)$ is near $F_0(\cdot)$ or $\hat{F}_Y(\cdot)$. 
These results are qualitatively similar to those for the one-sample, finite-dimensional example in \citet[\S4]{Kline2011}.

\begin{table}[htbp]
\centering
\caption{\label{tab:sim-RP}Bayesian test rejection probabilities, $\alpha=0.1$, $\numnornd{1000}$ replications.} 
\begin{tabular}{rrrr}
\toprule
&     &  \multicolumn{2}{c}{Comparison distribution} \\
\cmidrule{3-4}
$H_0$   &     
$n$ & 
$\UnifDist(0,1)$ & $Y$ \\
\midrule
    SD1 &    $10$ & 0.740 & 0.655 \\
    SD1 &    $40$ & 0.935 & 0.917 \\
    SD1 &   $100$ & 0.981 & 0.977 \\
    SD1 &$\infty$ & 1.000 & 1.000 \\[2pt]
non-SD1 &    $10$ & 0.000 & 0.005 \\
non-SD1 &    $40$ & 0.000 & 0.000 \\
non-SD1 &   $100$ & 0.000 & 0.000 \\
non-SD1 &$\infty$ & 0.000 & 0.000 \\
\bottomrule
\end{tabular}
\end{table}

\Cref{tab:sim-RP} shows rejection probabilities of the Bayesian test when $\theta(\cdot)=0(\cdot)$. 
This is the ``least favorable configuration'' for the null of SD1 (but not for non-SD1). 
The DGP has 
$X_i \iid \UnifDist(0,1)$ 
for $i=1,\ldots,n$. 
For one-sample testing, $F_0(\cdot)$ is the true (standard uniform) CDF of $X_i$. 
For two-sample testing, 
$Y_i \iid \UnifDist(0,1)$ for $i=1,\ldots,n$, identical to $X_i$. 
The hypotheses, methods, and notation are the same as for \cref{tab:sim-pval}. 
The entries for $n=\infty$ use \cref{prop:SD1-size,cor:nonSD1-RP}.

\Cref{tab:sim-RP} shows the same patterns as \cref{tab:sim-pval}.  
When $H_0$ is SD1, the Bayesian type I error rate is well above $\alpha$ even with $n=10$, with rejection probability increasing to $100\%$ as $n$ grows; consequently, size is also above $\alpha$. 
The opposite occurs when $H_0$ is non-SD1, which is not a subset of a half-space, with type I error rates of zero.%
\footnote{Although the type I error rate for non-SD1 is near zero with this DGP, the test's size is actually $\alpha$, which is attained when there is a single ``contact point'' with $F_X(r)=F_Y(r)$ and the inequalities are strict for all other $t \ne r$, thus reducing the test (asymptotically) to a single, scalar inequality.}

\end{document}